\documentclass[12pt]{amsart}
\usepackage{amsmath, amssymb}
\usepackage{multirow}

\newtheorem{theorem}{Theorem}[section]
\newtheorem{lemma}[theorem]{Lemma}

\newtheorem{conjecture}[theorem]{Conjecture}
\newtheorem{corollary}[theorem]{Corollary}
\theoremstyle{remark}
\newtheorem{remark}[theorem]{Remark}

\newcommand{\bK}{{\mathbb K}}
\newcommand{\bN}{{\mathbb N}}
\newcommand{\bQ}{{\mathbb Q}}
\newcommand{\bR}{{\mathbb R}}
\newcommand{\bZ}{{\mathbb Z}}

\newcommand{\cD}{{\mathcal D}}
\newcommand{\cN}{{\mathcal N}}

\def\disc{\operatorname{disc}}
\def\ord{\operatorname{ord}}

\numberwithin{equation}{section}

\begin{document}
\baselineskip=17pt

\title[Lang's conjecture and sharp height estimates for $y^{2}=x^{3}+ax$]{Lang's conjecture and sharp height estimates for the elliptic curves $y^{2}=x^{3}+ax$}

\author{Paul Voutier}
\address{}
\curraddr{London, UK}
\email{paul.voutier@gmail.com}

\author{Minoru Yabuta}
\address{}
\curraddr{Senri High School, 17-1, 2 chome, Takanodai, Suita, Osaka, 565-0861, Japan}
\email{yabutam@senri.osaka-c.ed.jp, rinri216@msf.biglobe.ne.jp}

\date{}

\dedicatory{}


\begin{abstract}
For elliptic curves given by the equation $E_{a}: y^{2}=x^{3}+ax$, we establish
the best-possible version of Lang's conjecture on the lower bound for the
canonical height of non-torsion rational points along with best-possible upper
and lower bounds for the difference between the canonical and logarithmic height.
\end{abstract}

\maketitle

\section{Introduction}

Heights, a measure of the arithmetic complexity of number-theoretic objects,
play a crucial role in the study of many diophantine problems (see \cite{BG}
for an excellent treatment of this subject). For points on elliptic curves,
there are two height functions of particular importance. Let $E$ be an elliptic
curve defined over a number field $\bK$ and denote by $E(\bK)$ the additive
group of all $\bK$-rational points on the curve $E$. For a point
$P \in E(\bK)$, we define the {\it canonical height} of $P$ by 
$$
\widehat{h}(P)=\frac{1}{2}\lim_{n \to \infty}\frac{h(2^{n}P)}{4^{n}},
$$
with $h(P)=h(x(P))$, where $h(P)$ and $h(x(P))$ are the {\it absolute logarithmic
heights} of $P$ and $x(P)$, respectively (see Sections~VIII.6,7 and 9 of
\cite{Silv2}). Also recall that for $\bQ$, $h(s/t)=\log\max \{|s|, |t|\}$
with $s/t$ in lowest terms is the absolute logarithmic height of $s/t$.

\subsection{Lower Bounds}

There are two types of lower bounds commonly used. We can look for
``Lehmer conjecture'' lower bounds where we fix the curve, $E$, and consider
lower bounds as the field of definition of the point, $P$, on $E$ varies
(see, for example, \cite{Masser}).

Alternatively, we can look for ``Lang conjecture'' lower bounds where we
fix the field of definition and vary the curve. It is such lower bounds that
we consider here for the elliptic curves, $E_{a}/\bQ$, given by the Weierstrass
equation $y^{2}=x^{3}+ax$ with $a \in \bZ$.

\begin{conjecture}[Lang's Conjecture]
Let $E/\bK$ be an elliptic curve with minimal discriminant $\cD_{E/\bK}$. There
exist constants $C_{1}>0$ and $C_{2}$, depending only on $[\bK:\bQ]$, such that
for all nontorsion points $P \in E(\bK)$ we have
$$
\widehat{h}(P) > C_{1} \log \left( \cN_{\bK/\bQ} \left( \cD_{E/\bK} \right) \right) + C_{2}.
$$
\end{conjecture}

See page~92 of \cite{Lang} along with the strengthened version in Conjecture~VIII.9.9
of \cite{Silv2}.

Such lower bounds have applications to counting the number of integral points
on elliptic curves (see \cite{HS}), problems involving elliptic divisibility
sequences \cite{Everest1, EIS, VY1}), \ldots

Silverman \cite{Silv1} showed that Lang's conjecture holds for any elliptic curve
with integral $j$-invariant over any number field (note that this includes our
curves, $E_{a}$, since their $j$-invariant is $1728$). Hindry and Silverman
\cite{HS} later proved an explicit version of Lang's conjecture whenever Szpiro's
ratio, $\sigma_{E/\bK}$, of $E/\bK$ is known. Hence Lang's conjecture follows
from Szpiro's conjecture (or the $ABC$ conjecture).

It can be shown that $\sigma_{E_{a}/\bQ}<4$, hence from Theorem~0.3 of \cite{HS},
$$
\widehat{h}(P) > 4 \cdot 10^{-33} \log \left( \cD_{E_{a}} \right).
$$

Subsequently, David \cite{David} and Petsche \cite{Petsche} improved Hindry
and Silverman's result. From Petsche's Theorem~2, for example, it follows that
$4 \cdot 10^{-33}$ above can be replaced by $10^{-21}$.

For $E_{a}/\bQ$ in the special case of $a=-n^{2}$ for a square-free integer
$n$, Bremner, Silverman and Tzanakis \cite[Proposition~2.1]{Bremner} proved a
much sharper result, namely,
$$
\widehat{h}(P) \geq \frac{1}{16}\log \left( 2n^{2} \right).
$$

In this paper, we provide a version for $E_{a}/\bQ$ for all non-zero integers
$a$ by examining $a$ which are fourth-power-free (i.e., global minimal
Weierstrass equations for all $E_{a}/\bQ$). Furthermore, our values of $C_{1}$
and $C_{2}$ are best-possible for $E_{a}/\bQ$.

\begin{theorem}
\label{thm:lang2}
Suppose $a$ is a fourth-power-free integer. Let $P \in E_{a}(\bQ)$ be a
nontorsion point. Then
$$
\widehat{h}(P) > \frac{1}{16}\log |a|
+ \left\{
	\begin{array}{ll}
		 (1/2)\log(2)    & \mbox{if $a>0$ and $a \equiv 1, 5, 7, 9, 13, 15 \bmod 16$} \\
		 (1/4)\log(2)    & \mbox{if $a>0$ and either $a \equiv 20,36 \bmod 64$} \\
		                 & \mbox{or $a \equiv 2,3,6,8,10,11,12,14 \bmod 16$} \\
		-(1/8)\log(2)    & \mbox{if $a>0$ and $a \equiv 4,52 \bmod 64$} \\
		 (9/16)\log(2)   & \mbox{if $a<0$ and $a \equiv 1, 5, 7, 9, 13, 15 \bmod 16$} \\
		 (5/16)\log(2)   & \mbox{if $a<0$ and  either $a \equiv 20,36 \bmod 64$} \\
		                 & \mbox{or $a \equiv 2,3,6,8,10,11,12,14 \bmod 16$} \\
		-(1/16)\log(2)   & \mbox{if $a<0$ and $a \equiv 4,52 \bmod 64$.}
	\end{array}
\right.
$$
\end{theorem}

\begin{remark}
For each of these cases, we have been able to find values of $a$ and points
$P \in E_{a}(\bQ)$ such that $\widehat{h}(P)$ is arbitrarily close to these
bounds, so they are best possible. We provide details in Section~\ref{sect:sharp}.
\end{remark}

Note that while the formulation of our result is not in terms of $\cD_{E_{a}}$,
it is equivalent to such a formulation since $\cD_{E_{a}}= \left| \Delta \left( E_{a} \right) \right|
=\left| -64a^{3} \right|$ for $a$ fourth-power-free. So we have the following
corollary.

\begin{corollary}
Suppose $a$ is a non-zero integer. If $P \in E_{a}(\bQ)$ is a nontorsion point,
then
$$
\widehat{h}(P) > \frac{1}{48}\log \left( \cD_{E_{a}} \right) - \frac{\log(2)}{4}.
$$
\end{corollary}

Our proof is based on the decomposition of the canonical height as the sum of
local height functions. To obtain our best-possible results, we require precise
bounds on the archimedean height on $E_{a}$ (Section~3, in particular,
Lemmas~\ref{lem:arch-a-neg} and \ref{lem:arch-a-pos}), along with a complete
analysis of the $p$-adic reduction of $E_{a}$ (Sections~4 and 5) and of the
denominators of $x(2P)$ (Section~6).

\subsection{Difference of Heights}

Due to the simple relationship between $h(P)$ and $P$, explicit bounds on
$(1/2)h(P)-\widehat{h}(P)$ are a key result for determining all points of bounded
canonical height on an elliptic curve. As a consequence, such bounds permit an
effective proof of the Mordell-Weil Theorem, and sharp bounds allow us to
determine Mordell-Weil bases of elliptic curves (see \cite[Chapter~X]{Silv2}).
In the same way, such bounds are also important for determining integral points
on elliptic curves.

Our proof of our lower bound for the canonical height also allows us to prove a
best possible upper bound on the difference between the canonical height and
the logarithmic height of points on $E_{a}(\bQ)$ and a very sharp lower bound.

In Example~2.2 of \cite{Silv5}, Silverman showed that
$$
-\frac{1}{4}\log |a| - 2.038
\leq \frac{1}{2}h(P)-\widehat{h}(P)
\leq \frac{1}{4}\log |a| + 2.252
$$
and that the coefficients on $\log|a|$ are best-possible.

Using a combination of Proposition~5.18(a) and Theorem~5.35(c) of \cite{S-Z}, one can obtain
$$
-\frac{1}{4}\log |a|-0.520 \leq \frac{1}{2}h(P)-\widehat{h}(P)
\leq \frac{1}{4}\log |a|+1.271.
$$

\begin{theorem}
\label{thm:hgt-diff}
Let $a$ be a non-zero integer. For all points $P \in E_{a}(\bQ)$, 
$$
-\frac{1}{4}\log |a| - \frac{1}{2\sqrt{|a|}}
< \frac{1}{2}h(P)-\widehat{h}(P)
< \frac{1}{4}\log |a|+ \frac{3}{8} \log(2).
$$

When $|a|$ is small, it is better to use the lower bound
$$
-\frac{1}{4}\log |a| - 0.16 < \frac{1}{2}h(P)-\widehat{h}(P).
$$
\end{theorem}

\begin{remark}
As for our lower bounds on the height, we have been able to find values of $a$
and points $P \in E_{a}(\bQ)$ such that the difference of the heights is
arbitrarily close to the upper and lower bounds stated here, so they are best
possible. We provide details in Section~\ref{sect:sharp}.
\end{remark}

\section{Notation}

For what follows in the remainder of this paper, we will require some standard
notation (see \cite[Chapter 3]{Silv2}, for example).

Let $\bK$ be a number field and let $E/\bK$ be an elliptic curve given by the
Weierstrass equation
$$
E: y^{2} + a_{1}xy + a_{3}y = x^{3} + a_{2}x^{2} + a_{4}x + a_{6},
$$
with $a_{1},\ldots,a_{6} \in \bK$.

Put
\begin{eqnarray*}
b_{2} & = & a_{1}^{2}+4a_{2}, \\
b_{4} & = & 2a_{4}+a_{1}a_{3}, \\
b_{6} & = & a_{3}^{2}+4a_{6}, \\
b_{8} & = & a_{1}^{2}a_{6}+4a_{2}a_{6}-a_{1}a_{3}a_{4}+a_{2}a_{3}^{2}-a_{4}^{2},
\end{eqnarray*}
then $E/\bK$ is also given by $y^{2}=4x^{3}+b_{2}x^{2}+2b_{4}x+b_{6}$.

Let $M_{\bK}$ be the set of valuations of $\bK$ and for each $v \in M_{\bK}$,
let $n_{v}$ be the local degree and let $\widehat{\lambda}_{v}(P):
E \left( \bK_{v} \right) \backslash \{ O \} \rightarrow \bR$ be the local
height function, where $\bK_{v}$ is the completion of $\bK$ at $v$. From
Theorem~VI.2.1 of \cite{Silv6}, we have the following decomposition of the
canonical height into local height functions
$$
\widehat{h}(P) = \sum_{v \in M_{\bK}} n_{v} \widehat{\lambda}_{v}(P).
$$

For $\bK=\bQ$, the non-archimedean valuations on $\bK$ can be identified with
the set of rational primes. For a non-archimedean valuation, $v$, we let
$q_{v}$ be the associated prime,
$$
v(x) = -\log |x|_{v} = \ord_{q_{v}}(x)\log \left( q_{v} \right)
$$
for $x \neq 0$ and $v(0)=+\infty$.

\begin{remark}
We refer the reader to \cite[Section~4]{Cremona} and \cite[Remark~VIII.9.2]{Silv2}
for notes about the various normalisations of both the canonical and local height
functions. In what follows, our local height functions, $\widehat{\lambda}_{v}(P)$,
are those that \cite{Cremona} denotes as $\lambda_{v}^{{\rm SilB}}(P)$, that is
as defined in Silverman's book \cite[Chapter~VI]{Silv6}. So as stated in (11) of
\cite{Cremona}, their $\lambda_{v} \left( P \right)$ equals
$2\widehat{\lambda}_{v}\left( P \right)+(1/6) \log \left| \Delta \left( E \right) \right|_{v}$
here.

Our canonical height also follows Silverman and is one-half that found in \cite{Cremona}
as well as one-half that returned from the height function, ellheight, in PARI.
\end{remark}

\section{Archimedean Estimates}

We will estimate the archimedean local height by using Tate's series (see
\cite{Tate} as well as the presentation in \cite{Silv3}). For any elliptic
curve, $E$, let
$$
t(P)=1/x(P) \hspace{3.0mm} \mbox{and} \hspace{3.0mm}
z(P) = 1-b_{4}t(P)^{2}-2b_{6}t(P)^{3}-b_{8}t(P)^{4},
$$
for a point $P=(x(P), y(P)) \in E(\bR)$. Then the archimedean local height of
$P \in E(\bR)$ is given by the series
\begin{equation}
\label{eq:arch-hgt}
\widehat{\lambda}_{\infty}(P)
= \frac{1}{2} \log |x(P)|
+ \frac{1}{8} \sum_{k=0}^{\infty} 4^{-k} \log |z(2^{k}P)|
- \frac{1}{12} \log \left| \Delta(E) \right|,
\end{equation}
provided $x \left( 2^{k}P \right) \neq 0$ for all $k \geq 0$.

\subsection{$a<0$}

\begin{lemma}
\label{lem:arch-a-neg}
Suppose $a \in \bR$ is negative and let $P=(x(P),y(P)) \in E_{a}(\bR)$ be a
point of infinite order.

\noindent
{\rm (a)} For $P \in E_{a}(\bR)$,
$$
\widehat{\lambda}_{\infty}(P) > \frac{1}{4} \log \left( x(P)^{2}-a \right)
- \frac{1}{12} \log \left| \Delta \left( E_{a} \right) \right|.
$$

\noindent
{\rm (b)} For $P \in E_{a}^{0}(\bR)$,
$$
-\frac{1}{3} \log(2)
< \left( \frac{1}{2} \log \max \left\{ 1, |x(P)| \right\}
-\frac{1}{12} \log \left| \Delta \left( E_{a} \right) \right| \right)
- \widehat{\lambda}_{\infty}(P)
< 0.
$$

\noindent
{\rm (c)} For $a \leq -2$ and $P \not\in E_{a}^{0}(\bR)$,
\begin{eqnarray*}
&      & -\frac{1}{4} \log |a| - \frac{1}{2\sqrt{|a|}} \\
&   <  & \left( \frac{1}{2} \log \max \left\{ 1, |x(P)| \right\}
-\frac{1}{12} \log \left| \Delta \left( E_{a} \right) \right| \right)
- \widehat{\lambda}_{\infty}(P) < -\frac{1}{4} \log (2) \nonumber.
\end{eqnarray*}

\noindent
{\rm (d)} For $a \leq -2$ and $P \not\in E_{a}^{0}(\bR)$,
\begin{equation}
\label{eq:a-neg-p-neg-arch2}
-\frac{1}{4} \log |a| - 0.16
< \left( \frac{1}{2} \log \max \left\{  1, |x(P)| \right\}
- \frac{1}{12} \log \left| \Delta \left( E_{a} \right) \right|
\right)
- \widehat{\lambda}_{\infty}(P).
\end{equation}
\end{lemma}

\begin{remark}
\label{rem:arch-a-neg}
The lower bound in part~(a) is approached as $x(P) \rightarrow +\infty$.

The lower bound in part~(b) is not sharp. It appears that the correct bound is
$-(1/4)\log(2)$ which is approached as $x(P)$ approaches $\sqrt{|a|}$. However,
the upper bound is attained as $x(P) \rightarrow +\infty$.

The lower bound in part~(c) is not sharp either. The correct bound appears to
be $-(1/4)\log|a|-0.22847\ldots |a|^{-1/2}$ which is attained at $x(P)=-1$ as
$a \rightarrow -\infty$. Note the coefficient of this $|a|^{-1/2}$ term
here equals the one mentioned in Remark~\ref{rem:arch-a-pos} for $a>0$.
The upper bound in part~(c) is attained at $x(P)=-\sqrt{|a|}$.

The lower bound in part~(d) is also not sharp. It appears that the correct
bound is $-(1/4)\log|a|-0.1310\ldots$, which is attained at $x(P)=-1$ for
$a=-2$. And as part~(c) demonstrates, this constant approaches $0$ as
$a \rightarrow -\infty$.
\end{remark}

\begin{proof}
If $x \left( 2^{k}P \right)=0$ for some non-negative integer $k$ (so that \eqref{eq:arch-hgt}
above does not converge), then $y\left( 2^{k}P \right)=0$ and $P$ is a torsion
point. But we are assuming that $P$ is not a torsion point, so $x\left( 2^{k}P \right) \neq 0$
for all such $k$.

For $E_{a}$, we have $a_{1}=a_{2}=a_{3}=a_{6}=0$, $a_{4}=a$, so $b_{2}=b_{6}=0$,
$b_{4}=2a$, $b_{8}=-a^{2}$,
$$
t(P)=1/x(P) \hspace{3.0mm} \mbox{and} \hspace{3.0mm} z(P)=(-at(P)^{2}+1)^{2}.
$$

For $a<0$, $E_{a}(\bR)$ has two components, and every point, $(x,y)$, in the
identity component $E_{a}^{0}(\bR)$ satisfies $x \geq \sqrt{|a|}$. From
Corollary~V.2.3.1 of \cite{Silv6}, $E_{a}(\bR) \cong (\bR/\bZ) \times (\bZ/2\bZ)$.
Therefore, $2P$, and $2^{k}P$ for all $k \geq 1$, is in $E_{a}^{0}(\bR)$.
Hence if $P \in E_{a}^{0}(\bR)$ is of infinite order, then
\begin{equation}
\label{eq:tate-params-bnds}
x(P) > \sqrt{|a|}, \hspace{3.0mm} 0 < t(P) < \frac{1}{\sqrt{|a|}},
\hspace{3.0mm}
1 < z(P) < 4
\end{equation}
(note that the point, $P$, with $x(P)=\sqrt{|a|}$ has order $2$).

(a) For every $P \in E_{a}(\bR)$ with $x \left( 2^{k}P \right) \neq 0$ for
all $k \geq 0$,
\begin{equation}
\label{eq:arch-hgt1}
\widehat{\lambda}_{\infty}(P)
= \frac{1}{4} \log \left( x(P)^{2}-a \right)
+ \frac{1}{8} \sum_{k=1}^{\infty} 4^{-k} \log \left| z \left( 2^{k} P \right) \right|
-\frac{1}{12} \log \left| \Delta \left( E_{a} \right) \right|,
\end{equation}
since $x(P)^{4}z(P)=\left( x(P)^{2}-a \right)^{2}$, noting that
$x(P)^{2}-a>0$.

Since $1 < z \left( 2^{k} P \right) < 4$ for $k \geq 1$,
\begin{equation}
\label{eq:d-bnds}
-\frac{1}{12} \log(2)
<
\left( \frac{1}{4} \log \left( x(P)^{2}-a \right) - \frac{1}{12} \log \left| \Delta \left( E_{a} \right) \right| \right)
- \widehat{\lambda}_{\infty}(P)
< 0,
\end{equation}
proving part~(a).

(b) From the bounds for $z(P)$ in \eqref{eq:tate-params-bnds},
$$
0 < \frac{1}{8} \sum_{k=0}^{\infty} 4^{-k} \log \left| z \left( 2^{k} P \right) \right|
< \frac{\log(2)}{3}
$$
and so part~(b) follows from \eqref{eq:arch-hgt}.

(c) For $P \not\in E_{a}^{0}(\bR)$, we have $-\sqrt{|a|} \leq x(P) \leq 0$
and so
$$
2 \leq \left( x(P)^{2}-a \right) /\max \left( 1, |x(P)| \right)^{2} \leq 1+|a|,
$$
with the min at $x(P)=-\sqrt{|a|}$ (such $P$ has order $2$) and the max at
$x(P)=-1$ (for $x(P)$ from $-1$ to $0$, it decreases from $1+|a|$ to $|a|$).
Hence
\begin{equation}
\label{eq:x4z-bnds}
\frac{1}{4} \log (2)
< \frac{1}{4} \log \left( x(P)^{2}-a \right)
- \frac{1}{2} \log \max \left\{ 1, |x(P)| \right\}
\leq \frac{1}{4} \log \left( 1+|a| \right),
\end{equation}
for any $P \not\in E_{a}^{0}(\bR)$ of infinite order.

From \eqref{eq:d-bnds}, along with the first inequality in \eqref{eq:x4z-bnds},
our upper bound in part~(c) follows.

From \eqref{eq:d-bnds}, the second inequality in \eqref{eq:x4z-bnds} and since
$\log \left( 1+|a| \right) \leq \log |a| + 1/|a|$, we obtain
\begin{eqnarray*}
&      & -\frac{1}{4} \log |a| - \frac{1}{4|a|} - \frac{1}{12}\log(2) \\
&   <  & \left( \frac{1}{2} \log \max \left\{ 1, |x(P)| \right\}
-\frac{1}{12} \log \left| \Delta \left( E_{a} \right) \right| \right)
- \widehat{\lambda}_{\infty}(P).
\end{eqnarray*}

But we can improve this lower bound.

If $x(P) \leq x_{1}=-\sqrt{2^{1/3}|a|/\left( \exp(1/\sqrt{|a|})|a|-2^{1/3}\right)}$,
then
$$
\left( x(P)^{2}-a \right) /\max \left( 1, |x(P)| \right)^{2}
\leq \exp \left(1/\sqrt{|a|} \right) |a|/2^{1/3}.
$$

Combining this inequality with our expression for $\widehat{\lambda}_{\infty}(P)$
in \eqref{eq:arch-hgt} and since $x(P)^{4}z(P)= \left( x(P)^{2}-a \right)^{2}$
and $1 < z \left( 2^{k}P \right) < 4$ for $k \geq 1$, we find for such $P$ that
\begin{eqnarray*}
\widehat{\lambda}_{\infty}(P)
& \leq & \frac{1}{4}\log|a| + \frac{1}{4\sqrt{|a|}}
+ \frac{1}{2} \log \max \left( 1, |x(P)| \right)
- \frac{1}{12} \log \left| \Delta \left( E_{a} \right) \right| \\
& & - \frac{1}{12}\log(2)
+ \frac{1}{8} \sum_{k=1}^{\infty} 4^{-k} \log \left| z \left( 2^{k}P \right) \right| \\
& \leq & \frac{1}{4}\log|a| + \frac{1}{4\sqrt{|a|}}
+ \frac{1}{2} \log \max \left( 1, |x(P)| \right)
- \frac{1}{12} \log \left| \Delta \left( E_{a} \right) \right|.
\end{eqnarray*}

Since $\left( \exp(1/\sqrt{|a|}) \right) |a| > |a| + \sqrt{|a|}$
and $\sqrt{|a|}> 2^{1/3}$ for $a \leq -2$, it follows that $x_{1} > -2^{1/6}$,
so the lower bound in part~(b) is satisfied for $a \leq -2$ and $P$ with
$x(P) \leq -2^{1/6}$.

So we now consider the remaining points with $-2^{1/6} < x(P) < 0$. From
\eqref{eq:arch-hgt1} and the second inequality in \eqref{eq:x4z-bnds}, it
suffices to bound from above
$$
\frac{1}{4} \log \left( 1 + 1/|a| \right)
+ \frac{1}{8} \sum_{k \geq 1} 4^{-k} \log \left| z \left( 2^{k}P \right) \right|
$$
and examine the sum
$$
Z(P) = \frac{1}{8} \sum_{k \geq 1} 4^{-k} \log \left| z \left( 2^{k}P \right) \right|.
$$

We can write $x(P)=-(1+\epsilon)$ where $-1 < \epsilon < 0.1225$ (note that
$2^{1/6}-1<0.1225$) and let $k_{0}$ be the smallest positive integer
such that $4^{k_{0}} \geq (1-\epsilon)\sqrt{|a|}$ (we will motivate this choice
of $k_{0}$ below).

We first bound the initial terms in the sum, $Z(P)$:
$$
Z_{0}(P)=\frac{1}{8}\sum_{k=1}^{k_{0}-1} 4^{-k} \log \left| z \left( 2^{k}P \right) \right|.
$$

Observe that
$$
x(2P) - \left( -\frac{(1-\epsilon)a}{4} \right)
= \frac{\epsilon^{2}a^{2}+ \left( \epsilon^{2}-3 \right)\left( 1+\epsilon \right)^{2}a+\left( 1+\epsilon \right)^{4}}
{4\left( 1+\epsilon \right)\left( -\left( 1+\epsilon \right)^{2}-a \right)}.
$$

For $a<0$, the numerator of the right-hand side will be positive as long as the
linear coefficient in $a$ is negative. That is, provided
$-\sqrt{3}<\epsilon<\sqrt{3}$. Similarly, the denominator is positive provided
that $-\sqrt{|a|}-1<\epsilon<\sqrt{|a|}-1$. Hence
$$
x(2P) > \left( \frac{(1-\epsilon)|a|}{4} \right),
$$
for all $-1 < \epsilon < 0.1225$.

If $x \geq \sqrt{|a|}$ and $a<0$, then $\left( x^{2}-a \right)/x>x$
and $\left( x^{2}-a \right)/\left( x^{2}+a \right)>1$. Hence
$x(2Q) > x(Q)/4$ for $Q \in E_{a}(\bR)$ with $x(Q) \geq \sqrt{|a|}$.

Therefore, $x \left( 2^{k}P \right) > 4^{-k}(1-\epsilon)|a|$ for all $k \geq 1$
and all $P \in E_{a}(\bR)$ with $x(P)=-(1+\epsilon)$ where $-1 < \epsilon < 0.1225$.
Hence, for all $1 \leq k<k_{0}$, $x \left( 2^{k}P \right) > 4^{-k}(1-\epsilon)|a|>\sqrt{|a|}$
-- and this is why we make our particular choice of $k_{0}$.

Since $a<0$ and $\log(1+x)<x$ for $x>0$, we have
$\log \left| z \left( 2^{k}P \right) \right| < 2|a|/x^{2}\left( 2^{k}P \right)$,
by the definition of $z\left( P \right)$. Since $x \left( 2^{k}P \right)>4^{-k}(1-\epsilon)|a|$,
$$
Z_{0}\left( P \right)
< \frac{1}{8}\sum_{k=1}^{k_{0}-1} 4^{-k} \frac{2|a|}{x^{2}\left( 2^{k}P \right)}
< \frac{4^{k_{0}-1}}{3(1-\epsilon)^{2}|a|}.
$$

Next, we bound the remaining terms. Since $z \left( 2^{k}P \right) \leq 4$,
$$
Z_{1}\left( P \right) = \frac{1}{8}\sum_{k \geq k_{0}} 4^{-k} \log \left| z \left( 2^{k}P \right) \right|
= \frac{1}{6}4^{-k_{0}} \log(4).
$$

Combining these two bounds and noting that we can write $4^{k_{0}}=c(1-\epsilon)\sqrt{|a|}$
with $1 \leq c < 4$,
\begin{equation}
\label{eq:z-bnd}
Z(P) < \left( \frac{c}{12} + \frac{\log(2)}{3c} \right) \frac{1}{(1-\epsilon)\sqrt{|a|}}
< \frac{0.3911}{(1-\epsilon)\sqrt{|a|}},
\end{equation}
since $c/12 + \log(2)/(3c)<0.3911$ for $1 \leq c < 4$ with the maximum value at $c=4$.

Hence, by \eqref{eq:z-bnd},
\begin{eqnarray*}
\widehat{\lambda}_{\infty}(P)
& < & \frac{1}{4}\log(1+|a|)
      + \frac{1}{2} \log \max \left( 1, |x(P)| \right)
      - \frac{1}{12} \log \left| \Delta \left( E_{a} \right) \right| \\
&   & + \frac{0.3911}{(1-\epsilon)\sqrt{|a|}}.
\end{eqnarray*}

To reduce the size of the constant in the $1/\sqrt{|a|}$ term here, we consider
two different ranges of $a$. 

For $-21 < -16/(1-\epsilon)^{2} \leq a \leq -2$, we have $k_{0}=1$, so $Z_{0}(P)=0$
for such $a$. Here
$$
\sqrt{|a|} \left( (1/4)\log(1+1/|a|) + Z_{1}(P) \right)
\leq \sqrt{|a|} \left( \frac{\log(1+1/|a|)}{4} + \frac{\log(4)}{24} \right)
< 0.318,
$$
where the maximum value is approached for $a$ near $-21$.

For $a \leq -21$, from \eqref{eq:z-bnd},
$$
(1/4)\log(1+1/|a|)+Z(P) < \frac{0.5}{\sqrt{|a|}}.
$$

(d) We use the same ranges for $a$ here as at the end of the proof of part~(c).

For $-21 < a \leq -2$,
$$
(1/4)\log(1+1/|a|) + Z_{1}(P) < (1/4)\log(1+1/|a|)+\log(4)/24 <0.1592,
$$
which obtained at $a=-2$.

Using the upper bound found in the proof of part~(c) for $a \leq -21$,
$1/(2\sqrt{|a|}) < 0.11$, and part~(d) follows.
\end{proof}

\subsection{$a>0$}

The following lemma will permit us to obtain sharp bounds for the archimedean
local height when $a>0$ in Lemma~\ref{lem:arch-a-pos} below.

\begin{lemma}
Let $a \in \bR$ be a positive real number and let
$$
E_{a}': y^{2}=x^{3}-3\sqrt{a}x^{2}+4ax-2a^{3/2}.
$$

For $P' \in E_{a}'(\bR)$, put
$$
z(P')=1-8ax(P')^{-2}+16a^{3/2}x(P')^{-3}-8a^{2}x(P')^{-4}.
$$

\noindent
{\rm (a)}
Suppose that $x(P')=c\sqrt{a}$ where $c>1$. Then
\begin{equation}
\label{eq:x4z-lb1}
\frac{1}{8} \log \left( \frac{x(P')^{4}z(P')}{x(P)^{4}} \right)
> \frac{1}{2c} - \frac{3}{4c^{2}}.
\end{equation}

\noindent
{\rm (b)}
Let $k$ be a positive integer and suppose that $x(P') \geq c\sqrt{a}$ where
$c > 4^{k-1}$. Then
\begin{equation}
\label{eq:logz2kp-lb1}
\log \left( z \left( 2^{k}P' \right) \right)
> -\frac{128 \cdot 4^{2(k-1)}}{c^{2}}.
\end{equation}

\noindent
{\rm (c)}
Suppose that $x(P')=(1+\epsilon)\sqrt{a}$ where $\epsilon>0$. Then
\begin{equation}
\label{eq:x4z-lb2}
\frac{1}{8} \log \left( x(P')^{4}z(P') \right)
> \frac{1}{4} \log |a| + \frac{\epsilon}{2} - \frac{5\epsilon^{2}}{4}.
\end{equation}

\noindent
{\rm (d)}
Let $k$ be a positive integer and suppose that $x(P')=(1+\epsilon)\sqrt{a}$,
where $0<\epsilon<4^{-(k-1)}$. Then
\begin{equation}
\label{eq:logz2kp-lb2}
\log \left( z \left( 2^{k}P' \right) \right)
> -128 \cdot 4^{2(k-1)} \epsilon^{2}.
\end{equation}
\end{lemma}

\begin{proof}
We start with some bounds for the logarithm function that will be used
throughout this section.

Using the Taylor series expansion of $\log(1+x)$, we have
$\log(1+x) \geq x-x^{2}/2+x^{3}/3-x^{4}/4$ for $0 \leq x <1$. Since the
derivative of the lower bound, $-x^{3}+x^{2}-x+1$, has a negative leading
coefficient and its largest real root is at $x=1$ where $\log(1+x)=\log(2)$
while the value of the lower bound is $5/12$, this lower bound also applies
for $x \geq 1$. Hence, for $x \geq 0$,
\begin{equation}
\label{eq:log-plus-lb1}
\log(1+x) \geq x-x^{2}/2+x^{3}/3-x^{4}/4.
\end{equation}

Similarly, for $x \geq 0$,
\begin{equation}
\label{eq:log-plus-lb2}
\log(1+x) \geq x-x^{2}/2.
\end{equation}

By its Taylor series expansion, for $0 \leq x \leq 0.6$,
\begin{equation}
\label{eq:log-neg-lb}
\log (1-x) \geq -x-x^{2}.
\end{equation}

(a) Substituting $x(P')=c\sqrt{a}$ into the expression for $x(P')^{4}z(P')$,
and recalling that $x(P')=x(P)+\sqrt{a}$, we obtain
\begin{equation}
\label{eq:x4z1}
\frac{x(P')^{4}z(P')}{x(P)^{4}}
=\frac{c^{4}-8c^{2}+16c-8}{(c-1)^{4}}=1+\frac{4c^{3}-14c^{2}+20c-9}{(c-1)^{4}}.
\end{equation}

The numerator of the last term in \eqref{eq:x4z1} has positive leading
coefficient and one real root, at $c=0.7835\ldots$, hence that last term is
positive for $c>0.79$.

Applying \eqref{eq:log-plus-lb1} to our expression for $x(P')^{4}z(P')/x(P)^{4}$
in \eqref{eq:x4z1}, we find that $\log \left( x(P')^{4}z(P')/x(P)^{4} \right)$
is bounded below by a rational function in $c$ whose numerator, $f_{n}(c)$, is
a polynomial of degree $15$ in $c$ with $48$ as its leading coefficient and
whose denominator is $f_{d}(c)=12(c-1)^{16}$. Using Maple, we find that the
polynomial $c^{2}f_{n}(c)-c^{2}(4/c-6/c^{2})f_{d}(c)$ has $208$ as its leading
coefficient and that its largest positive root is at $4.6564\ldots$, and so
this polynomial is positive for such values of $c$. Since $f_{d}(c)>0$ for such
$c$, part~(a) holds for such $c$.

For $1<c<4.6564\ldots$, we observe that the left-hand side of \eqref{eq:x4z-lb1}
exceeds the right-hand side near $c=1$, since $x(P')^{4}z(P')/x(P)^{4}$ grows
arbitrarily large as $c$ approaches $1$ from the right, while the right-hand
side of \eqref{eq:x4z-lb1} is negative. The derivative of the left-hand side is
$$
-\frac{c^{3}-4c^{2}+8c-4}{2(c^{4}-8c^{2}+16c-8)(c-1)},
$$
which is negative for $c>1$, since the numerator and denominator have no roots
for such $c$.

The derivative of the left-hand-side minus the right-hand side of \eqref{eq:x4z-lb1}
is
$$
-\frac{13c^{4}-52c^{3}+96c^{2}-80c+24}{2(c^{4}-8c^{2}+16c-8)(c-1)c^{3}}.
$$
Once again, this is negative for $c>1$, since the numerator and denominator
have no roots for such $c$. Therefore, since \eqref{eq:x4z-lb1} holds for $c$
close to $1$ and for $c>4.66$, it must hold for all $1<c<4.66$ too, concluding
the proof of part~(a).

(b) We start by showing that if $x(P')=c\sqrt{a}$ where $c>1$, then
\begin{equation}
\label{eq:x-lb}
x(2P') > \frac{c}{4}\sqrt{a}.
\end{equation}

We can write
$$
x(2P')= \sqrt{a} \frac{c^{4}-8c^{2}+16c-8}{4\left( c^{3}-3c^{2}+4c-2 \right)}.
$$

Now
$$
c^{4}-8c^{2}+16c-8 = 4 \left( c^{3}-3c^{2}+4c-2 \right) \left( \frac{c}{4} + 0.455 \right)
$$
simplifies to $1.18c^{3}-6.54c^{2}+10.72c-4.36=0$, which only has one real root
at $0.6066\ldots$. Furthermore, $c^{3}-3c^{2}+4c-2>0$ for $c>1$, so in fact,
for all $c>1$, a stronger inequality than \eqref{eq:x-lb} holds,
$x(2P') > (c/4+0.455)\sqrt{a}$.

Furthermore, since the right-hand side is monotonically increasing in $c$,
for any $Q'$ with $x(Q') \geq c\sqrt{a}$ and $c>1$, we have
$x(2Q')>(c/4+0.455)\sqrt{a}$.

We now prove \eqref{eq:logz2kp-lb1}, proceeding by induction.

We start with $k=1$.

Again we suppose that $x \left( P' \right) = c \sqrt{a}$ (rather than
greater than or equal to $c\sqrt{a}$). Substituting the exact expression for
$x(2P')$ into the expression for $z(P')$, we find that $z(2P')=1-z_{2}(P')$,
where
$$
z_{2}(P') = \frac{128c^{4} \left( c-1 \right)^{2} \left( c^{2}-2c+2 \right)^{2}\left( c-2 \right)^{4}}
{\left( c^{4}-8c^{2}+16c-8 \right)^{4}}.
$$

Since $z_{2}(P') \geq 0$, the numerator of $z_{2}(P')$ minus $3/5$ times its
denominator is a polynomial in $c$ of degree $16$ with $-3/5$ as its leading
coefficient and its largest real root at $0.95456\ldots$, we have $0 \leq z_{2}(P')<0.6$
and so
$$
\log(z(2P')) \geq -z_{2}(P')-z_{2}^{2}(P'),
$$
by \eqref{eq:log-neg-lb}.

Using Maple, we can write $-z_{2}(P')-z_{2}^{2}(P')-(-128/c^{2})$ as a
rational function of $c$ where the numerator is of degree $31$ with $1792$ as
its leading coefficient and the denominator is
$\left( c^{4}-8c^{2}+16c-8 \right)^{8}c^{2}$.
The largest real root of the numerator occurs at $0.8988\ldots$, the largest
root of the denominator occurs at $0.78315\ldots$ and so
\begin{equation}
\label{eq:logz2Pa}
\log \left( z \left( 2P' \right) \right)
\geq -\frac{128}{c^{2}},
\end{equation}
for all $c \geq 1$.

We now show that the result holds for $x(P') \geq c\sqrt{a}$.

First observe that since
$$
\frac{d}{dx(P')}z(P')=\frac{16a\left(x^{2}-3\sqrt{a}x+2a \right)}{x^{5}}
= \frac{16a\left(x-2\sqrt{a} \right) \left( x-\sqrt{a} \right)}{x^{5}},
$$
we find that $z(P')$ has a maximum at $x(P')=\sqrt{a}$
where $z(P')=1$, a minimum at $x(P')=2\sqrt{a}$ where $z(P')=1/2$ and $z(P')$
asymptotically approaches $1$ from below for $x(P')>2\sqrt{a}$. Hence
\begin{equation}
\label{eq:zp-bnds}
1/2 \leq z(P') \leq 1.
\end{equation}

For $1 \leq c < 13.58$, $-128/c^{2} < -\log(2)$. Thus, for $c$ in this range,
(b) holds. Using Maple, we find that for $c>5.62$, $z(2P')$ is monotonically
increasing. Hence (b) holds for all larger values of $c$ as well, completing
the proof for $k=1$.

For $k>1$, applying \eqref{eq:x-lb} repeatedly, we find that
$x \left( 2^{k}P' \right) > 4^{-k}c\sqrt{a}$ for $c>4^{k-1}$. Hence part~(b)
follows from \eqref{eq:logz2Pa}.

(c) Note that $x(P')^{4}z(P')=x(P')^{4}-8ax(P')^{2}+16a^{3/2}x(P')-8a^{2}$.
So if $x(P')=(1+\epsilon)\sqrt{a}$ for $\epsilon>0$, then
\begin{equation}
\label{eq:x4z2}
x(P')^{4}z(P')
= \left( 1 + 4\epsilon - 2\epsilon^{2} + 4\epsilon^{3}+\epsilon^{4} \right) a^{2}.
\end{equation}

For $0<\epsilon<2$, $0<\epsilon_{1}=4\epsilon - 2\epsilon^{2}<4\epsilon - 2\epsilon^{2}+4\epsilon^{3}+\epsilon^{4}$
and hence, by \eqref{eq:log-plus-lb2},
\begin{eqnarray*}
\log \left( x(P')^{4}z(P')/a^{2} \right) & > & \log \left( 1+\epsilon_{1} \right) > \epsilon_{1} - \epsilon_{1}^{2}/2 \\
& = & 4\epsilon-10\epsilon^{2}+8\epsilon^{3}-2\epsilon^{4} \\
& > & 4\epsilon-10\epsilon^{2},
\end{eqnarray*}

For $\epsilon \geq 2$, $4\epsilon-10\epsilon^{2}<0$ while
$1 + 4\epsilon - 2\epsilon^{2} + 4\epsilon^{3}+\epsilon^{4}>1$, so part~(c) follows.

(d) We proceed by induction and similarly to part~(b).

We will also show that for $k \geq 1$ and all $0<\epsilon<4^{-(k-1)}$,
\begin{equation}
\label{eq:x2kp}
x \left( 2^{k}P' \right) > \left( 4^{-k}/\epsilon \right) \sqrt{a}.
\end{equation}

With $x(P')=(1+\epsilon)\sqrt{a}$ for $\epsilon>0$, we have
\begin{equation}
\label{eq:x2p}
x(2P') = \frac{1+4\epsilon-2\epsilon^{2}+4\epsilon^{3}+\epsilon^{4}}
{4\epsilon \left( 1+\epsilon^{2} \right)}\sqrt{a}.
\end{equation}

Since $\left( 1+4\epsilon-2\epsilon^{2}+4\epsilon^{3}+\epsilon^{4} \right)
-\left( 1+\epsilon^{2} \right) = 4\epsilon-3\epsilon^{2}+4\epsilon^{3}+\epsilon^{4}$
has no positive real roots, for $\epsilon>0$,
$$
x(2P') > \frac{1}{4\epsilon} \sqrt{a}.
$$

Using Maple, we substitute \eqref{eq:x2p} into the definition of $z(2P')$ and
find that $z(2P')$ is the quotient of a monic polynomial in $\epsilon$ of degree
$16$ with integer coefficients divided by
$\left( 1+4\epsilon-2\epsilon^{2}+4\epsilon^{3}+\epsilon^{4} \right)^{4}$.
We denote the numerator and denominator as $z_{2,n}(\epsilon)$ and
$z_{2,d}(\epsilon)$, respectively.

Taking the series expansion of this expression for $z(2P')$, we can write
it as $1-128\epsilon^{2}+2048\epsilon^{3}-\cdots$. Using Maple again, we
find that the polynomial $z_{2,n}(\epsilon)-z_{2,d}(\epsilon) \left( 1-128\epsilon^{2}+1100\epsilon^{3} \right)$,
which is of degree $19$ in $\epsilon$ and whose leading coefficient is $-1000$,
has a triple root at $\epsilon=0$, a root at $0.07475\ldots$ and its other
three real roots are negative. Since $z(2P')>1-128\epsilon^{2}+1100\epsilon^{3}$
for $\epsilon=0.02$, it follows that
$$
z(2P')>1-128\epsilon^{2}+1100\epsilon^{3}
$$
for all $0 < \epsilon < 0.07475$.

Applying \eqref{eq:log-neg-lb} with $x=128\epsilon^{2}-1100\epsilon^{3}$,
we obtain a polynomial that we can show (e.g., using Maple) is larger than
$-128\epsilon^{2}$ in the desired range and thus
$$
\log \left( z(2P') \right) > - 128\epsilon^{2},
$$
for $0<\epsilon<0.074$. For $\epsilon \geq 0.074$, we have $-128\epsilon^2<-0.7<-\log(2)$,
and so, from \eqref{eq:zp-bnds}, part~(d) holds for $k=1$ and any $\epsilon>0$.

Now suppose that $k \geq 2$ and that $x(P')=(1+\epsilon)\sqrt{a}$. Put
$c=4^{-(k-1)}/\epsilon$. By our inductive hypothesis, $x \left( 2^{k-1}P' \right) \geq c\sqrt{a}$,
since our assumption on $\epsilon$ ensures that $\epsilon<4^{-(k-2)}$.
From \eqref{eq:x-lb}, it follows that \eqref{eq:x2kp} holds. Applying
\eqref{eq:logz2kp-lb1} with $k$ there set to $1$ and our value of
$c$ here completes the proof of \eqref{eq:logz2kp-lb2}.
\end{proof}

\begin{lemma}
\label{lem:arch-a-pos}
Let $a \in \bR$ be a positive real number and $P=(x(P),y(P)) \in E_{a}(\bR)$
be a point of infinite order.

\noindent
{\rm (a)}
\begin{equation}
\label{eq:archim-hgt-a-pos2}
\widehat{\lambda}_{\infty}(P)
> \frac{1}{4} \log(a) - \frac{1}{12} \log \left| \Delta \left( E_{a} \right) \right|.
\end{equation}

\noindent
{\rm (b)} For $a \geq 2$,
\begin{eqnarray*}
&   & - \frac{1}{4}\log(a) - \frac{1}{2\sqrt{a}}  \\
& < & \left( \frac{1}{2} \log \max \left\{ 1, \left| x(P) \right| \right\}
             -\frac{1}{12} \log \left| \Delta \left( E_{a} \right) \right|
      \right)
      - \widehat{\lambda}_{\infty}(P) < 0. \nonumber
\end{eqnarray*}

\noindent
{\rm (c)} For $a \geq 3$,
$$
- \frac{1}{4}\log(a) - 0.16
< \left( \frac{1}{2} \log \max \left\{ 1, \left| x(P) \right| \right\}
-\frac{1}{12} \log \left| \Delta \left( E_{a} \right) \right|
\right)
- \widehat{\lambda}_{\infty}(P).
$$
\end{lemma}

\begin{remark}
\label{rem:arch-a-pos}
The lower bound in part~(a) is approached as $x(P)$ approaches $0$.

The correct lower bound in part~(b) appears to be
$-(1/4)\log(a)-0.26033\ldots a^{-1/2}$ which is attained at $x(P)=1$ with
$a=2$. In fact, $-(1/4)\log(a)-0.22847\ldots a^{-1/2} -Ca^{-1}$ for an absolute
constant $C>0$ appears to hold. Observe that the coefficient of this $a^{-1/2}$
term here equals the one mentioned in Remark~\ref{rem:arch-a-neg} for $a<0$
as $a \rightarrow -\infty$. The upper bound in part~(b) is attained at
$x(P) \rightarrow +\infty$.

The lower bound in part~(c) is also not sharp. It appears that the correct
bound is $-(1/4)\log|a|-0.14922\ldots$, which is attained at $x(P)=1$ for
$a=3$. As part~(b) demonstrates, this constant approaches $0$ as
$a \rightarrow +\infty$.
\end{remark}

\begin{proof}
(a) For $a>0$, $E_{a}(\bR)$ has only one component and it includes
$(0,0)$ which causes a problem since we require $x(2^{k}P)$ to be bounded away
from $0$ to ensure that Tate's series converges. To get around this, we use an
idea of Silverman's (see page~340 of \cite{Silv3}) and translate the curve to
the right using $x'=x+\sqrt{a}$, noting that $\widehat{\lambda}_{\infty}$ is
fixed under such translations (i.e., if $E'$ is the translated curve with
$\widehat{\lambda}_{\infty}'(P)$ as the archimedean local height function, then
$\widehat{\lambda}_{\infty}'(P')=\widehat{\lambda}_{\infty}(P)$). In this way,
we obtain the equation
$$
E_{a}': y^{2}=x^{3}-3\sqrt{a}x^{2}+4ax-2a^{3/2}
$$
and every point, $(x,y)$, in $E_{a}'(\bR)$ satisfies $x \geq \sqrt{a}$.

For $E_{a}'$, we have $b_{2}=-12\sqrt{a}$, $b_{4}=8a$, $b_{6}=-8a^{3/2}$ and
$b_{8}=8a^{2}$. Hence
$$
t(P')=1/x(P') \hspace{3.0mm} \mbox{and} \hspace{3.0mm}
z(P')=1-8at(P')^{2}+16a^{3/2}t(P')^{3}-8a^{2}t(P')^{4}
$$
for any $P' \in E_{a}'(\bR)$.

From \eqref{eq:zp-bnds}, for every $P' \in E_{a}'(\bR)$, 
\begin{eqnarray}
\label{eq:init-bnd}
\widehat{\lambda}_{\infty}(P')
& = & \frac{1}{8} \log \left| x(P')^{4}z(P') \right|
+ \frac{1}{8} \sum_{k=1}^{\infty} 4^{-k} \log \left| z \left( 2^{k}P' \right) \right|
-\frac{1}{12} \log \left| \Delta \left( E_{a} \right) \right| \nonumber \\
& \geq & \frac{1}{8} \log \left| x(P')^{4}z(P') \right| - \frac{\log(2)}{24}
-\frac{1}{12} \log \left| \Delta \left( E_{a} \right) \right|,
\end{eqnarray}
since $-\log(2) \leq \log \left| z \left( 2^{k}P' \right) \right| \leq 0$
and $\Delta \left( E_{a} \right) = \Delta \left( E_{a}' \right)$.

In fact, the term $-\log(2)/24$ is unnecessary. We use a more careful analysis
here to show that.

Write $x(P')=(1+\epsilon)\sqrt{a}$ where $\epsilon>0$.

First, suppose that $x(P') \geq 1.07\sqrt{a}$. From \eqref{eq:x4z2} and our
observation that $x(P')^{4}z(P')$ increases as $\epsilon$ does, it follows that
$x(P')^{4}z(P')>1.27a^{2}$.

Hence $(1/8) \log \left( x(P')^{4}z(P') \right)-(1/24)\log(2)
>(1/4)\log(a)+\log(1.27)/8-\log(2)/24>(1/4)\log(a)$. So, from \eqref{eq:init-bnd},
$$
\widehat{\lambda}_{\infty}(P')
> \frac{1}{4}\log(a)-\frac{1}{12} \log \left| \Delta \left( E_{a} \right) \right|.
$$

Now suppose that $\epsilon<0.07$.

For $N \geq 1$ and $0<\epsilon<(2/7) \cdot 4^{-N}$, then
\begin{eqnarray*}
\widehat{\lambda}_{\infty}(P')
& > & \frac{1}{4}\log(a)+\frac{\epsilon}{2}-\frac{5\epsilon^{2}}{4}
-\left( \sum_{k=N+1}^{\infty} \frac{\log(2)}{8 \cdot 4^{k}}
+ \sum_{k=1}^{N} \frac{128 \cdot 4^{2(k-1)} \epsilon^{2}}{8 \cdot 4^{k}} \right)
-\frac{1}{12} \log \left| \Delta \left( E_{a} \right) \right| \\
& = & \frac{1}{4}\log(a)+\frac{\epsilon}{2}-\frac{4^{N+2}-1}{12}\epsilon^{2}
-\frac{\log(2)}{6 \cdot 4^{N+1}}
-\frac{1}{12} \log \left| \Delta \left( E_{a} \right) \right|,
\end{eqnarray*}
by \eqref{eq:x4z-lb2} and \eqref{eq:logz2kp-lb2}.

Define
$$
f_{N}(x) = -\frac{4^{N+2}-1}{12}x^{2} + \frac{x}{2} - \frac{\log(2)}{6 \cdot 4^{N+1}}.
$$

$f_{N}(x)$ has two distinct real roots, $\alpha_{N}$ and $\beta_{N}$ with $\alpha_{N}<\beta_{N}$.
So if $\alpha_{N} < x < \beta_{N}$, then $f_{N}(x)>0$. Put $t_{N}=(2/7) \cdot 4^{-N}$.
For $N \geq 1$,
$$
f_{N} \left( t_{N} \right) = \frac{4^{N} (40-49\log(2))+8}{1176 \cdot 4^{2N}} > 0
$$
and
$$
f_{N} \left( t_{N+1} \right) = \frac{4^{N} (68-98\log(2))+1}{2352 \cdot 4^{2N}} > 0.
$$

Therefore, $\alpha_{N} < t_{N+1} < t_{N} < \beta_{N}$.
Hence if $t_{N+1} \leq x \leq t_{N}$, then $f_{N}(x)>0$.

Thus these intervals overlap and cover all $0<\epsilon \leq 0.07$, we
conclude that
$$
\widehat{\lambda}_{\infty}(P) = \widehat{\lambda}_{\infty}(P')
> \frac{1}{4} \log(a) - \frac{1}{12} \log \left| \Delta \left( E_{a} \right) \right|
$$
holds for all $P'$ with $x(P')>\sqrt{a}$.

(b) Note that
$$
x(P')^{4}z(P') =
x(P)^{4} + 4\sqrt{a}x(P)^{3} - 2ax(P)^{2} + 4a^{3/2}x(P) + a^{2}.
$$

Now
$$
\frac{d}{dx} \frac{x^{4} + 4\sqrt{a}x^{3} - 2ax^{2} + 4a^{3/2}x + a^{2}}{x^{4}}
=-4\frac{\sqrt{a}x^{3} -ax^{2}+3a^{3/2}x+a^{2}}{x^{5}},
$$
and its numerator only has one real root, near $-0.2955\ldots \sqrt{a}$.

Hence
$$
\frac{x^{4} + 4\sqrt{a}x^{3} - 2ax^{2} + 4a^{3/2}x + a^{2}}{x^{4}}
$$
is monotonically decreasing towards $1$ for $x>0$ and is
$1+4\sqrt{a}-2a+4a^{3/2}+a^{2}$ at $x=1$.

Similarly,
$$
x^{4} + 4\sqrt{a}x^{3} - 2ax^{2} + 4a^{3/2}x + a^{2}
$$
is monotonically increasing for $x \geq 0$ with the value $a^{2}$
at $x=0$ and $1+4\sqrt{a}-2a+4a^{3/2}+a^{2}$ at $x=1$.

Therefore,
$$
1 \leq \frac{x(P')^{4}z(P')}{\max \left\{ 1, \left| x(P) \right| \right\}^{4}}
\leq 1+4\sqrt{a}-2a+4a^{3/2}+a^{2}.
$$

Applying this to our expression for $\widehat{\lambda}_{\infty}(P')$ in
\eqref{eq:arch-hgt}, we obtain
\begin{eqnarray}
\label{eq:arch-at-1}
&      & -\frac{1}{8} \log \left( 1+4\sqrt{a}-2a+4a^{3/2}+a^{2} \right) \\
& \leq & \left(
\frac{1}{2} \log \max \left\{ 1, \left| x(P) \right| \right\}
+ \frac{1}{8} \sum_{k=1}^{\infty} 4^{-k} \log \left| z \left( 2^{k}P' \right) \right|
-\frac{1}{12} \log \left| \Delta \left( E_{a} \right) \right| \right) \nonumber \\
&      & - \widehat{\lambda}_{\infty}(P) \leq 0. \nonumber
\end{eqnarray}

Since $1/2 \leq z \left( 2^{k}P' \right) \leq 1$ and
$(1/8)\log \left( 1+4\sqrt{a}-2a+4a^{3/2}+a^{2} \right)
< (1/2)\log \left( 1+\sqrt{a} \right)
= (1/4)\log (a) + (1/2) \log \left( 1+1/\sqrt{a} \right)
< (1/4)\log (a) + (1/2)/\sqrt{a}$,
it follows that
\begin{eqnarray*}
&   & -\frac{1}{4} \log(a) - \frac{1}{2\sqrt{a}} \\
& < & \left( \frac{1}{2} \log \max \left\{ 1, \left| x(P) \right| \right\}
             -\frac{1}{12} \log \left| \Delta \left( E_{a} \right) \right|
      \right)
      - \widehat{\lambda}_{\infty}(P)
\leq \frac{1}{24} \log(2).
\end{eqnarray*}

Our lower bound in part~(b) follows.

But as in part~(a), we can show that the term $\log(2)/24$ in the upper bound
is unnecessary.

We will show that
\begin{equation}
\label{eq:arch-lb-a-pos}
\hspace*{5.0mm}
0 < \frac{1}{8} \log \left( \frac{x(P')^{4}z(P')}{\max \left\{ 1, \left| x(P) \right| \right\}^{4}} \right)
    + \frac{1}{8} \sum_{k=1}^{\infty} 4^{-k} \log \left| z \left( 2^{k}P' \right) \right|.
\end{equation}

First consider $0 \leq x(P) \leq 1$. We showed earlier in the proof of part~(b)
that $(1/8)\log \left| x(P')^{4}z(P') \right| \geq (1/4)\log(a)$ for such $x(P)$.
For $a \geq 2$, $(1/4)\log(a)>\log(2)/24$ and so \eqref{eq:arch-lb-a-pos} holds.

Now consider $x(P) \geq 1$ and write $x(P')=c\sqrt{a}$ where $c>1$.

For $N \geq 1$ and $c>4^{N-1}$, then
\begin{eqnarray*}
&   & \frac{1}{8} \log \left( \frac{x(P')^{4}z(P')}{\max \left\{ 1, \left| x(P) \right| \right\}^{4}} \right)
      + \frac{1}{8} \sum_{k=1}^{\infty} 4^{-k} \log \left| z \left( 2^{k}P' \right) \right| \\
& > & \frac{1}{2c}-\frac{3}{4c^{2}}
-\left( \sum_{k=N+1}^{\infty} \frac{\log(2)}{8 \cdot 4^{k}}
+ \sum_{k=1}^{N} \frac{128 \cdot 4^{2(k-1)}}{8 \cdot 4^{k} c^{2}} \right) \\
& = & \frac{\epsilon}{2}-\frac{4^{N+2}-25}{12}\epsilon^{2}
-\frac{\log(2)}{6 \cdot 4^{N+1}},
\end{eqnarray*}
by \eqref{eq:x4z-lb1} and \eqref{eq:logz2kp-lb1}.

Letting $\epsilon=1/c$, we can now define $f_{N}(x)$ as in the proof of
part~(a) and proceed in the same way as there.

(c) For $a \geq 3$, we have
$(1/8)\log \left( 1+4\sqrt{a}-2a+4a^{3/2}+a^{2} \right)
< (1/4)\log(a)+0.16$. Part~(c) follows from this, \eqref{eq:arch-at-1} and the
fact that $z \left( 2^{k}P' \right) \leq 1$.
\end{proof}

\section{Non-archimedean Estimates for $q_{v}$ odd}

\begin{lemma}
\label{lem:nonarch-v-odd}
Let $v$ be a non-archimedean valuation on $\bQ$ associated with an odd
prime number, $q_{v}$, and let $a$ be an integer such that $q_{v}^{4} \nmid a$.

\noindent
{\rm (a)} The Kodaira symbols and Tamagawa indices of $E_{a}$ at $v$ are as in
Table~$\ref{table:kodaira-odd}$.

\begin{table}[h]
\centering
\begin{tabular}{|c|c|c|}
\hline
$a$                  & {\rm Kodaira symbol} & {\rm Tamagawa index} \\ \hline
$\ord_{q_{v}}(a)=0$  & $I_{0}$              & $1$                               \\ \hline
$\ord_{q_{v}}(a)=1$  & $III$                & $2$                               \\ \hline
$\ord_{q_{v}}(a)=2$, &                      & \\
{\rm Legendre symbol} $\left( \left( -a/q_{v}^{2} \right)/q_{v} \right)=1$  & $I_{0}^{*}$ & $4$          \\ \hline
$\ord_{q_{v}}(a)=2$, &                      & \\
{\rm Legendre symbol} $\left( \left( -a/q_{v}^{2} \right)/q_{v} \right)=-1$ & $I_{0}^{*}$ & $2$          \\ \hline
$\ord_{q_{v}}(a)=3$  & $III^{*}$            & $2$                               \\ \hline
\end{tabular}
\caption{$E_{a}$ reduction information for $q_{v}$ odd}
\label{table:kodaira-odd}
\end{table}

\noindent
{\rm (b)} For any $P \in E_{a} \left( \bQ_{v} \right)$, $2P$ is always non-singular
and if $2P \neq O$, then
\begin{equation}
\label{eq:padic-2p-hgt}
\widehat{\lambda}_{v}(2P) = \frac{1}{2} \log \max\{ 1, \left| x(2P) \right|_{v} \}
- \frac{\log \left| \Delta \left( E_{a} \right) \right|_{v}}{12}.
\end{equation}

\noindent
{\rm (c)} For any $P \in E_{a}\left( \bQ_{v} \right) \backslash \{ O \}$,
\begin{eqnarray*}
\widehat{\lambda}_{v}(P) & = &\frac{1}{2} \log \max\{ 1, \left| x(P) \right|_{v} \} 
- \frac{ \log \left| \Delta \left( E_{a} \right) \right|_{v}}{12} \\
& & - \left\{
	\begin{array}{ll}
		(1/4)\log \left( q_{v} \right) & \mbox{if $q_{v}||a$ and $\ord_{q_{v}}(x(P))>0$}, \\
		(1/2)\log \left( q_{v} \right) & \mbox{if $q_{v}^{2}||a$ and $\ord_{q_{v}}(x(P))>0$}, \\
		(3/4)\log \left( q_{v} \right) & \mbox{if $q_{v}^{3}||a$ and $\ord_{q_{v}}(x(P))>0$}, \\
		0 & \mbox{otherwise.}
	\end{array}
\right.
\end{eqnarray*}

\noindent
{\rm (d)} $P \in E_{a} \left( \bQ_{v} \right)$ is singular if and only if
$\ord_{q_{v}}(x(P)), \ord_{q_{v}}(y(P))>0$.
\end{lemma}

\begin{proof}
(a) We use Tate's algorithm with $K=\bQ_{v}$ to obtain the reduction information
below (using the steps and notation in Silverman's presentation of Tate's
algorithm in Section~IV.9 of \cite{Silv6}).

\noindent
$\bullet$ Step~1. This step applies when $\ord_{q_{v}}\left( \Delta \left( E_{a} \right) \right) =0$,
so the Kodaira symbol is $I_{0}$ at $v$ when $\ord_{q_{v}}(a)=0$.

\noindent
$\bullet$ Step~2. We have $\ord_{q_{v}}\left( \Delta \left( E_{a} \right) \right)>0$.
The singular point, $P=(x(P), y(P))$, is already at $(0,0)$ since
$\ord_{q_{v}}(2y(P)), \ord_{q_{v}}(a)>0$ implies that $\ord_{q_{v}}(x(P))>0$ too,
so no change of variables is needed. Therefore, $b_{2}=0$ and hence
$\ord_{q_{v}} \left( b_{2} \right)>0$. Thus Step~2 does not apply.

\noindent
$\bullet$ Step~3. Since $a_{6}=0$ and hence $\ord_{q_{v}} \left( a_{6} \right) \geq 2$,
this step does not apply.

\noindent
$\bullet$ Step~4. If $\ord_{q_{v}}(a)=1$, then $\ord_{q_{v}} \left( b_{8} \right)=2<3$,
since $b_{8}=-a^{2}$. So Step~4 applies since $\ord_{q_{v}} \left( a_{6} \right) \geq 2$.
Therefore, the Kodaira symbol is $III$ at $v$ when $\ord_{q_{v}}(a)=1$.

\noindent
$\bullet$ Step~5. Step~5 does not apply since $b_{6}=0$ and hence
$\ord_{q_{v}} \left( b_{6} \right) \geq 3$.

\noindent
$\bullet$ Step~6. Since $\ord_{q_{v}}(a)=1$ is treated by Step~4, we now have
$\ord_{q_{v}}(a) \geq 2$.

Step~6 applies since $a_{1}=a_{2}=a_{3}=a_{6}=0$, and since $\ord_{q_{v}}(a) \geq 2$,
no change of coordinates is necessary. Thus $P(T)=T^{3}+(a/q_{v}^{2})T$ and
$\disc(P)=-4a^{3}/q_{v}^{6}$. If $\ord_{q_{v}}(a)=2$, then $\disc(P)$ is not
divisible by $q_{v}$. Hence the Kodaira symbol is $I_{0}^{*}$ at $v$ when
$\ord_{q_{v}}(a)=2$.

Note that the Tamagawa index, $c_{q_{v}}$, is $4$ if $-a/q_{v}^{2}$ is a 
quadratic residue modulo $q_{v}$ and is $2$ otherwise.

\noindent
$\bullet$ Step~7. We now have $\ord_{q_{v}}(a) \geq 3$. This implies that
$P(T)$ has a triple root and hence Step~7 does not apply.

\noindent
$\bullet$ Step~8. Step~8 does not apply since $a_{3}=a_{6}=0$ implies that
$Y^{2}+a_{3,2}Y-a_{6,4}=Y^{2}$ does not have two distinct roots, where, as in
Section~IV.9 of \cite{Silv6}, we use $a_{i,j}$ to denote $q_{v}^{-j}a_{i}$.

\noindent
$\bullet$ Step~9. This applies since $\ord_{q_{v}}(a)=3<4$. Thus the Kodaira
symbol is $III^{*}$ at $v$ when $\ord_{q_{v}}(a)=3$.

Since $q_{v}^{4} \nmid a$, this completes the proof of part~(a).

(b) From the characterisation of $E_{a} \left( \bQ_{v} \right)/E_{a}^{0}\left( \bQ_{v} \right)$
in Table~4.1 of \cite[Chapter IV]{Silv6}, we see that
$2E_{a}\left( \bQ_{v} \right) \subseteq E_{a}^{0} \left( \bQ_{v} \right)$ for
these Kodaira symbols and thus $2P$ is always non-singular. Hence
we have \eqref{eq:padic-2p-hgt} from Theorem~4.1 of \cite[Chapter VI]{Silv6}.

(c) This follows from using Table~\ref{table:kodaira-odd} here and Table~2 of
\cite{Cremona} (recall the difference between our local heights and those in
\cite{Cremona}).

(d) Finally, we determine when $P$, considered as an element of $E_{a} \left( \bQ_{v} \right)$,
is singular for $q_{v}|a$. We require $\ord_{q_{v}} \left( 3x(P)^2+a \right)>0$
and $\ord_{q_{v}} \left( y(P) \right)>0$. The second inequality along with
$q_{v}|a$ implies that $\ord_{q_{v}} \left( x(P) \right)>0$.
\end{proof}

\section{Non-archimedean Estimates for $q_{v}=2$}

\begin{lemma}
\label{lem:nonarch-v-2}
Let $a$ be an integer and suppose that $16 \nmid a$.

{\rm (a)} The Kodaira symbols and Tamagawa indices of $E_{a}$ at $2$ are as in
Table~$\ref{table:kodaira-2}$.

\noindent
{\rm (b)} For any $P \in E_{a} \left( \bQ_{2} \right)$, $4P$ is always
non-singular and if $2P \neq O$, then
\begin{eqnarray}
\label{eq:2adic-2p-hgt}
\widehat{\lambda}_{2}(2P) & = &\frac{1}{2} \log \max\{ 1, \left| x(2P) \right|_{2} \} 
- \frac{ \log \left| \Delta \left( E_{a} \right) \right|_{2}}{12}
\nonumber \\
& & - \left\{
	\begin{array}{ll}
		(1/2) \log(2) & \mbox{if $a \equiv 4, 52 \bmod 64$ and $\ord_{2}(x(2P))>0$} \\
		0 & \mbox{otherwise.}
	\end{array}
\right.
\end{eqnarray}

\noindent
{\rm (c)} For any $P \in E_{a}\left( \bQ_{2} \right) \backslash \{ O \}$,
\begin{eqnarray*}
\widehat{\lambda}_{2}(P) & = &\frac{1}{2} \log \max\{ 1, \left| x(P) \right|_{2} \} 
- \frac{ \log \left| \Delta \left( E_{a} \right) \right|_{2}}{12} \\
& & - \left\{
	\begin{array}{ll}
		(1/4)\log(2) & \mbox{if $a \equiv 2,3 \bmod 4$ and $\ord_{2}(x(P)+a)>0$}, \\
		(1/2)\log(2) & \mbox{if $a \equiv 12,20,36,44 \bmod 64$ and $\ord_{2}(x(P))>0$}, \\
		             & \mbox{or if $a \equiv 4, 28, 52, 60 \bmod 64$ and $\ord_{2}(x(P))>1$}, \\
		(3/4)\log(2) & \mbox{if $a \equiv 0 \bmod 8$ and $\ord_{2}(x(P))>0$} \\
		             & \mbox{or if $a \equiv 28, 60 \bmod 64$ and $\ord_{2}(x(P))=1$}, \\
		(7/8)\log(2) & \mbox{if $a \equiv 4, 52 \bmod 64$ and $\ord_{2}(x(P))=1$}, \\
		           0 & \mbox{otherwise.}
	\end{array}
\right.
\end{eqnarray*}

\noindent
{\rm (d)} $P \in E_{a} \left( \bQ_{2} \right)$ is singular if and only if
$\ord_{2}(x(P)+a)>0$.
\end{lemma}

\begin{table}[h]
\centering
\begin{tabular}{|c|c|c|}
\hline
$a$                       & {\rm Kodaira symbol} & {\rm Tamagawa index} \\ \hline
$a \equiv     1 \bmod  4$ & $II$                 & $1$                  \\ \hline
$a \equiv     3 \bmod  4$ & $III$                & $2$                  \\ \hline
$a \equiv     2 \bmod  4$ & $III$                & $2$                  \\ \hline
$a \equiv    12 \bmod 32$ & $I_{2}^{*}$          & $2$                  \\ \hline
$a \equiv    28 \bmod 32$ & $I_{2}^{*}$          & $4$                  \\ \hline
$a \equiv 20,36 \bmod 64$ & $I_{3}^{*}$          & $2$                  \\ \hline
$a \equiv  4,52 \bmod 64$ & $I_{3}^{*}$          & $4$                  \\ \hline
$a \equiv     0 \bmod  8$ & $III^{*}$            & $2$                  \\ \hline
\end{tabular}
\caption{$E_{a}$ reduction information for $q_{v}=2$}
\label{table:kodaira-2}
\end{table}

\begin{proof}
(a) For $q_{v}=2$, Tate's algorithm provides the reduction information below.

\noindent
$\bullet$ Step~1. Since $\Delta \left( E_{a} \right) = -64a^{3}$ is always
even, Step~1 never applies.

\noindent
$\bullet$ Step~2. If $a$ is even, then the singular point is already at $(0,0)$
so no translation is required. In this case, $b_{2}=a_{1}^{2}+4a_{2}=0$,
so Step~2 does not apply.

If $a$ is odd, then the singular point is at $(1,0)$ so we must use the change
of variables $x=x'+1$ and we have $a_{1}=a_{3}=0$, $a_{2}=3$,
$a_{4}=a+3$ and $a_{6}=a+1$. In this case, $b_{2}=12 \equiv 0 \bmod 2$ and
again Step~2 does not apply.

\noindent
$\bullet$ Step~3. If $a$ is even, then Step~3 does not apply, since $a_{6}=0$.

If $a$ is odd, then Step~3 only applies if $a \equiv 1 \bmod 4$, since
$a_{6}=a+1$. So the Kodaira symbol is $II$ at $v$ when $a \equiv 1 \bmod 4$.

\noindent
$\bullet$ Step~4. If $a$ is even, then $b_{8}=-a^{2}$. Hence if $\ord_{2}(a)=1$,
then the Kodaira symbol is $III$ at $v$.

If $a \equiv 3 \bmod 4$, then $b_{8}=3+6a-a^{2} \equiv 12 \bmod 16$. So the
Kodaira symbol is $III$ at $v$ in this case too.

The only remaining case is $a \equiv 0 \bmod 4$.

\noindent
$\bullet$ Step~5. Since $a$ is even, $b_{6}=0$, so this step does not apply.

\noindent
$\bullet$ Step~6. Here $P(T)=T^{3}+a_{4,2}T$ and $\disc(P)=-4a^{3}/64$,
recalling that, as in Section~IV.9 of \cite{Silv6}, we use $a_{i,j}$ to denote
$q_{v}^{-j}a_{i}$.

We know that $a \equiv 0 \bmod 4$. Hence $\disc(P) \equiv 0 \bmod 2$ and
so Step~6 does not apply.

\noindent
$\bullet$ Step~7. In this step, we must have $a \equiv 4 \bmod 8$ (otherwise,
$P(T)$ has a triple root). In this case, $P(T)$ has a double root at $T=1$ and
we apply the change of variables $x=x'+2$ to obtain a new equation for the
elliptic curve with $a_{2}=6$, $a_{3}=0$, $a_{4}=a+12$ and $a_{6}=2a+8$. Since
$Y^{2}+a_{3,2}Y-a_{6,4}=Y^{2}-a_{6,4}$ never has distinct roots, the
Kodaira symbol can never be $I_{1}^{*}$.

If $a \equiv 12 \bmod 16$, then $a_{6} \equiv 0 \bmod 32$ and hence
$Y^{2}+a_{3,2}Y-a_{6,4}=Y^{2}$ has its double root at $Y=0$.

In this case, we consider the polynomial $a_{2,1}X^{2}+a_{4,3}X+a_{6,5}$.
Its discriminant is $a_{4,3}^{2}-4a_{2,1}a_{6,5}=(a+12)^{2}/64-12(2a+8)/32$.
Writing $a=16a_{c}+12$, this discriminant is $4a_{c}^{2}-3 \not\equiv 0 \bmod 2$.
Hence $a_{2,1}X^{2}+a_{4,3}X+a_{6,5}$ has distinct roots and the Kodaira symbol
is $I_{2}^{*}$ at $v$ when $a \equiv 12 \bmod 16$.

Furthermore, $a_{2,1}X^{2}+a_{4,3}X+a_{6,5}$ is reducible (and so $c_{2}=4$)
if and only if $a_{6,5}=0$ (i.e., $a_{6} \equiv 0 \bmod 64$). This is equivalent
to $a \equiv 28 \bmod 32$. Hence $c_{2}=2$ if and only if $a \equiv 12 \bmod 32$.

If $a \equiv 4 \bmod 16$, then $a_{6} \equiv 16 \bmod 32$ and hence
$Y^{2}+a_{3,2}Y-a_{6,4}=Y^{2}-1$ has its double root at $Y=1$.

In this case, we apply the change of variables $y=y'+4$ to obtain a new
equation for the elliptic curve with $a_{2}=6$, $a_{3}=8$, $a_{4}=a+12$ and
$a_{6}=2a-8$. Since $Y^{2}+a_{3,2}Y-a_{6,4}=Y^{2}$ never has distinct roots,
the Kodaira symbol can never be $I_{1}^{*}$.

We now consider $a_{2,1}X^{2}+a_{4,3}X+a_{6,5}=3X^{2}+a_{6,5}$, which also
never has distinct roots, so the Kodaira symbol can never be $I_{2}^{*}$.

If $a \equiv 4 \bmod 32$, then $a_{4} \equiv 16 \bmod 32$ and
$a_{6} \equiv 0 \bmod 64$. Thus $a_{4,3} \equiv a_{6,5} \equiv 0 \bmod 2$,
so the double root of $a_{2,1}X^{2}+a_{4,3}X+a_{6,5}=a_{2,1}X^{2}$ is at $X=0$
and we consider $Y^{2}+a_{3,3}Y-a_{6,6}=Y^{2}+Y+a_{6,6}$. Its discriminant is
$1-4a_{6,6} \not\equiv 0 \bmod 2$ and hence it always has distinct roots, so
its Kodaira symbol is $I_{3}^{*}$ at $v$.

If $a \equiv 4 \bmod 64$, then $a_{6,6} \equiv 0 \bmod 128$. Hence
$Y^{2}+a_{3,3}Y-a_{6,6}=Y^{2}+Y$ and $c_{2}=4$.

If $a \equiv 36 \bmod 64$, then $a_{6,6} \equiv 64 \bmod 128$. Hence
$Y^{2}+a_{3,3}Y-a_{6,6}=Y^{2}+Y+1$ and $c_{2}=2$.

If $a \equiv 20 \bmod 32$, then the double root of $a_{2,1}X^{2}+a_{4,3}X+a_{6,5}$
is at $X=1$, so we need to apply a change of variables again. Here, we apply
the change of variables $x'=x''+4$ to obtain a new equation for the elliptic
curve with $a_{2}=18$, $a_{3}=8$, $a_{4}=a+108$ and $a_{6}=6a+200$. Thus
$a_{2,1}X^{2}+a_{4,3}X+a_{6,5}=a_{2,1}X^{2}$ has a double root at $X=0$ and we
consider $Y^{2}+a_{3,3}Y-a_{6,6}=Y^{2}+Y+a_{6,6}$. Its discriminant is
$1-4a_{6,6} \not\equiv 0 \bmod 2$ and hence it always has distinct roots,
so its Kodaira symbol is $I_{3}^{*}$ at $v$.

If $a \equiv 20 \bmod 64$, then $a_{6,6} \equiv 64 \bmod 128$. Hence
$Y^{2}+a_{3,3}Y-a_{6,6}=Y^{2}+Y+1$ and $c_{2}=2$.

If $a \equiv 52 \bmod 64$, then $a_{6,6} \equiv 0 \bmod 128$. Hence
$Y^{2}+a_{3,3}Y-a_{6,6}=Y^{2}+Y$ and $c_{2}=4$.

\noindent
$\bullet$ Step~8. $P(T)=T^{3}+a_{4,2}T$ has a triple root (and at $T=0$)
only if $a \equiv 0 \bmod 8$. However, $Y^{2}+a_{3,2}Y-a_{6,4}=Y^{2}$ does not
have distinct roots so Step~8 does not apply.

\noindent
$\bullet$ Step~9. Since $Y^{2}+a_{3,2}Y-a_{6,4}=Y^{2}$, this step does apply,
provided that $2^{4} \nmid a$ (which we have assumed). Therefore, the Kodaira
symbol is $III^{*}$ and $c_{2}=2$ at $v$ when $\ord_{2}(a)=3$.

Since $16 \nmid a$, this concludes the proof of part~(a) of the lemma.

(b) According to Table~4.1 of \cite{Silv6}, we see that $2P$ is non-singular
(as an element of $E_{a} \left( \bQ_{2} \right)$) unless the Kodaira symbol of
$E_{a}$ at $2$ is $I_{3}^{*}$ and $c_{2}=4$ (these conditions imply that the
component group is isomorphic to $\bZ/(4\bZ)$), which only happens for
$a \equiv 4, 52 \bmod 64$.

So for $a \not\equiv 4, 52 \bmod 64$ or $a \equiv 4,52 \bmod 64$ and $2P$
non-singular, we can apply Theorem~4.1 of \cite[Chapter VI]{Silv6} again.

For $a \equiv 4,52 \bmod 64$ and $2P$ is singular, we appeal to the case of
Kodaira symbol $I_{m}^{*}$, $m$ odd and $c_{v}=2$ or $4$ in the proof of
Proposition~6 of \cite{Cremona}. Our $2P$ here must be of order $2$ in
$E_{a} \left( \bQ_{2} \right)/E_{a}^{0} \left( \bQ_{2} \right)$ (since
$4P \in E_{a}^{0} \left( \bQ_{2} \right))$ and hence it must equal $P_{1}$ in
their proof. They calculate that their $\lambda_{v} \left( P_{1} \right) = -\log \left( q_{v} \right)/n_{v}$.
Since $n_{v}=1$ and $q_{v}=2$ here, their $\lambda_{v} \left( P_{1} \right)=-\log(2)$.
Recall that $\lambda_{v} \left( P_{1} \right)$ equals
$2\widehat{\lambda}_{v}\left( P_{1} \right)+(1/6) \log \left| \Delta \left( E_{a} \right) \right|_{2}$.
Since $2P$ is singular, from part~(d) we have $\left| x(2P) \right|_{2}<1$, so
$$
\widehat{\lambda}_{v}(2P) = \frac{1}{2} \log \max\{ 1, \left| x(2P) \right|_{2} \}
- \frac{ \log \left| \Delta \left( E_{a} \right) \right|_{2}}{12}
- \frac{1}{2}\log(2).
$$

(c) Except for the cases when the Kodaira symbol is $I_{m}^{*}$ and $c_{2}=4$,
this follows from using Table~\ref{table:kodaira-2} here and Table 2 of
\cite{Cremona}.

When $I_{m}^{*}$ and $c_{2}=4$, then Table~2 of \cite{Cremona} has two possible
values. Working through this case of the proof of Proposition~6 in \cite{Cremona}
(note that the authors of \cite{Cremona} consider $m$ even and $m$ odd separately
and divide each into subcases depending on whether $c_{v}=2$ or $c_{v}=4$, so
working through the proof for our cases is easy to do). In this way, we find
that the points that are labelled $P_{1}$ in \cite{Cremona} are those with
$\ord_{2}(x(P))>1$ and those labelled $P_{2}$ in \cite{Cremona} are those with
$\ord_{2}(x(P))=1$. The values of $\widehat{\lambda}_{2}(P)$ in these cases are
stated in the proof of Proposition~6 of \cite{Cremona}.

Again, recall the difference between our local heights and those in
\cite{Cremona}.

(d) So it remains to determine when $2P$ is singular as an element of
$E_{a} \left( \bQ_{2} \right)$. We require $\ord_{2} \left( 3x(P)^2+a \right)>0$
and $\ord_{2} \left( 2y(P) \right)>0$.

If $a$ is even, then $\ord_{2} \left( 3x(P)^2+a \right)>0$ is equivalent to
$\ord_{2} \left( x(P) \right)>0$. Since $\ord_{2} \left( x(P) \right)>0$
implies $\ord_{2} \left( y(P) \right)>0$, $\ord_{2}(x(P))>0$ is a necessary
and sufficient condition.

If $a$ is odd, then $\ord_{2} \left( 3x(P)^2+a \right)>0$ is equivalent to
$\ord_{2} \left( x(P) \right)=0$, which implies that
$\ord_{2} \left( y(P) \right)>0$, so $\ord_{2}(x(P))=0$ is a necessary and
sufficient condition here.

Both cases can be covered by the condition $\ord_{2}(x(P)+a)>0$.
\end{proof}

\section{$B_{2}$ estimates}

In addition to allowing us to estimate $x(2P)$ in the proof of Theorem~\ref{thm:lang2},
the following lemma will also be useful for our study of primitive divisors
(see \cite{VY1}).

\begin{lemma}
\label{lem:2.4}
Let $a$ be a fourth-power-free integer and let $P \in E_{a}(\bQ)$ be a point of
infinite order.

\noindent
{\rm (a)} Let $x(P)=uv^{2}$ with $u \in \bZ$ square-free and $v \in \bQ$.
If $n$ is even, then $x(nP)$ is a rational square. If $n$ is odd, then
$x(nP)=uw^{2}$ for some $w \in \bQ$.

\noindent
{\rm (b)} Writing $x(nP)=A_{n}/B_{n}$ in lowest terms with $A_{n} \in \bZ$
and $B_{n} \in \bN$, we have
$$
\ord_{2} \left( B_{2} \right)
\geq \left\{
	\begin{array}{ll}
		 4 & \mbox{if $a \equiv 1, 5, 7, 9, 13, 15 \bmod 16$} \\
		 2 & \mbox{if $a \equiv 2,3,6,8,10,11,12,14 \bmod 16$ or $a \equiv 20,36 \bmod 64$} \\
		 0 & \mbox{if $a \equiv 4,52 \bmod 64$}.
	\end{array}
\right.
$$

\noindent
{\rm (c)} Furthermore, if $a \not\equiv 4,52 \bmod 64$ or if $\ord_{2}(x(P)) \neq 1$,
then $\ord_{2} \left( B_{2} \right) \geq \ord_{2} \left( B_{1} \right) + 2$.
\end{lemma}


\begin{proof}
(a) Let $\bQ^{\ast}$ be the multiplicative group of non-zero rational numbers,
and let $\bQ^{\ast 2}$ denote the subgroup of squares of elements of
$\bQ^{\ast}$. We define a map $\alpha$ from $E_{a}(\bQ)$ to
$\bQ^{\ast}/\bQ^{\ast 2}$ by 
\begin{eqnarray*}
& & \alpha(O)=1, \;  \alpha((0,0))=a, \\
& & \alpha((x,y))=s \quad \text{if $x=st^{2}$ with $s$ square-free,}
\end{eqnarray*}
where $O$ is the zero element in $E_{a}(\bQ)$. Then $\alpha$ is
homomorphism (see p.85 \cite{Silv-Tate}). Let $x(P)=uv^{2}$ with $u \in \bZ$
square-free and $v \in \bQ$. Then 
\begin{eqnarray*}
\alpha(2P) & = & \alpha(P+P)=\alpha(P)^{2}=1, \\
\alpha(3P) & = & \alpha(2P+P)=\alpha(2P)\alpha(P)=u.
\end{eqnarray*}

Using induction shows that if $n$ is even, then $\alpha(nP)=1$, and if $n$ is
odd then $\alpha(nP)=u$. Therefore, if $n$ is even, then $x(nP)$ is a rational
square, and if $n$ is odd, then $x(nP)=uw^{2}$ for some $w \in \bQ$.

(b) From the arguments on pages~92--93 of \cite{Silv-Tate}, we can write
$P= \left( b_{1}M^{2}/e^{2}, b_{1}MN/e^{3} \right)$ in lowest terms,
where $a=b_{1}b_{2}$ with $\gcd \left( b_{1}, e \right)
=\gcd \left( b_{2}, M \right)=\gcd(e,M)=\gcd(M,N)=\gcd \left( e,N \right)=1$.

Suppose a prime $p$ divides both $b_{1}$ and $N$. Since $P \in E_{a}(\bQ)$,
we can write $N^{2} = b_{1}M^{4} + b_{2}e^{4}$. If $p^{2}|b_{1}$, then as
$\gcd \left( e,N \right)=1$, it follows that $p^{2}|b_{2}$. Therefore $p^{4}$
divides $a$, which contradicts the assumption that $a$ is fourth-power-free.
Hence $\ord_{p} \left( \gcd \left( b_{1}, N \right) \right) \leq 1$.

Since for any $Q=(x,y)$, by the duplication formula, we have 
$$
x(2Q)=\frac{(x^{2}-a)^{2}}{4y^{2}}=\frac{\left( 2x^{3}-y^{2} \right)^{2}}{4x^{2}y^{2}}.
$$

Applying these with the expression we stated above for $P$, we obtain
\begin{equation}
\label{eq:x2p-b2}
x(2P)=\frac{\left( b_{1}M^{4}-b_{2}e^{4} \right)^{2}}{4M^{2}N^{2}e^{2}}
=\frac{\left( 2b_{1}M^{4}-N^{2} \right)^{2}}{4M^{2}N^{2}e^{2}}.
\end{equation}

\noindent
\underline{Case 1:} If $N$ is odd, then $2b_{1}M^{4}-N^{2}$ is odd.

\noindent
Case 1-a: If $e$ is even, then $M$ is odd, since $\gcd(e,M)=1$.
Hence $\ord_{2} \left( B_{2} \right) = \ord_{2} \left( e^{2} \right)+2$. So
$\ord_{2} \left( B_{2} \right) = \ord_{2} \left( B_{1} \right) + 2 \geq 4$.

In this case, $b_{1}$ is odd, since $\gcd \left( b_{1}, e \right)=1$ and we
find that $a$ can take any non-zero value modulo $16$.

\noindent
Case 1-b: If $e$ is odd and $M$ is odd, then
$\ord_{2} \left( B_{2} \right) = \ord_{2} \left( B_{1} \right)+2=2$.

In this case, we find that $a$ can take any even non-zero value modulo $16$.

\noindent
Case 1-c: If $e$ is odd and $M$ is even, then
$\ord_{2} \left( B_{2} \right) = \ord_{2} \left( e^{2} \right)+2+2\ord_{2}(M)$.

So $\ord_{2} \left( B_{2} \right) \geq \ord_{2} \left( B_{1} \right)+4=4$.

In this case, we find that $a$ can take any non-zero value modulo $16$.

\noindent
\underline{Case 2:} Next suppose that $2\|N$, then $e$ and $M$ are both odd.

\noindent
Case 2-a: Suppose that $2 \nmid b_{1}$.  Then $b_{1}M^{4}$ is odd.
Hence, $4 \| \left( 2b_{1}M^{4}-N^{2} \right)^{2}$,
but $2^{4}\| \left( 4M^{2}N^{2}e^{2} \right)$,
since $\gcd(e,N)=\gcd(M,N)=1$. So $\ord_{2} \left( B_{2} \right) =
\ord_{2} \left( e^{2} \right)+2 = \ord_{2} \left( B_{1} \right)+2=2$.

We have $N^{2} \equiv b_{1}M^{4}+b_{2}e^{4} \bmod 16$. Since $2 \parallel N$,
$N^{2} \equiv 4 \bmod 16$. Since $e$ and $M$ are both odd,
$M^{4} \equiv e^{4} \equiv 1 \bmod 16$. Therefore $b_{1}+b_{2} \equiv 4 \bmod 16$.
Examining each of the possibilities for $b_{1}$ and $b_{2}$, we find that
$a=b_{1}b_{2} \equiv 3,11 \bmod 16$.

\noindent
Case 2-b: Suppose that $2 \| b_{1}$. In this case, we have
$N^{2} \equiv 4 \bmod 32$, $e$ and $M$ are both odd and
$b_{1}M^{4} \equiv b_{1} \bmod 32$ (since $b_{1}$ is even). Since
$N^{2} = b_{1}M^{4}+b_{2}e^{4} \equiv 4 \bmod 32$, it similarly follows that
$b_{2}e^{4} \equiv b_{2} \bmod 32$. Hence $b_{2}e^{4} - b_{1}M^{4} \equiv 0 \bmod 8$.
In this case, $\left( b_{1}M^{4}-b_{2}e^{4} \right)^{2}$ is divisible by $64$,
whereas $16 \| 4M^{2}N^{2}e^{2}$. Therefore $B_{2}$ is odd.

Substituting $b_{1}=4b_{11}+2$ and $b_{2} \equiv 4-b_{1} \bmod 32$ into
$a=b_{1}b_{2}$, we find that $a \equiv 48b_{11}^{2}+4 \bmod 64$. This is
either $4$ or $52 \bmod 64$ depending on whether $b_{11}$ is even or odd,
respectively.

\noindent
\underline{Case 3:} If $4|N$, then $4 \nmid b_{1}$ and both $e$ and $M$ are odd.

\noindent
Case 3-a: Suppose that $b_{1}$ is odd, then $b_{1}M^{4}$ is odd. Hence,
$\ord_{2} \left( \left( 2b_{1}M^{4}-N^{2} \right)^{2} \right)=2$, but
$\ord_{2}\left( 4M^{2}N^{2}e^{2} \right)=\ord_{2} \left( N^{2} \right)+2 \geq 6$,
since $\gcd(e,N)=1$. So, in this case,
$\ord_{2} \left( B_{2} \right) = \ord_{2} \left( B_{1} \right)+4 \geq 4$.

Examining each of the possibilities for $b_{1}$ and $b_{2}$, we find that
$b=b_{1}b_{2} \equiv 7, 15 \bmod 16$.

\noindent
Case 3-b: Suppose that $b_{1} \equiv 2 \bmod 4$, then
$\ord_{2} \left( \left( 2b_{1}M^{4}-N^{2} \right)^{2} \right) = 4$, but
$\ord_{2}\left( 4M^{2}N^{2}e^{2} \right)=\ord_{2} \left( N^{2} \right)+2 \geq 6$,
since $\gcd(e,N)=1$. So, in this case too,
$\ord_{2} \left( B_{2} \right) \geq \ord_{2} \left( e^{2} \right)+2
=\ord_{2} \left( B_{1} \right)+2=2$. 

Examining each of the possibilities for $b_{1}$ and $b_{2}$, we find that
$b=b_{1}b_{2} \equiv 12 \bmod 16$.

We now combine these results.

From case~(1-a), we find that $\ord_{2} \left( B_{2} \right) \geq 4$ always
holds. However, from case~(1-b), if $a$ is even, then $\ord_{2} \left( B_{2} \right)=2$
is possible. From case~(2-a), if $a \equiv 3, 11 \bmod 16$, then $\ord_{2} \left( B_{2} \right)=2$
is also possible. From case~(2-b), if $a \equiv 4,52 \bmod 64$, then
$\ord_{2} \left( B_{2} \right)=0$ is possible.

(c) For general $x(P)$ and $a \not\equiv 4,52 \bmod 64$, this follows from our
proof of part~(b).

Suppose first that $\ord_{2}(x(P)) \equiv 0 \bmod 2$. Recall
from the proof of part~(b) that we can write $P= \left( b_{1}M^{2}/e^{2},
b_{1}MN/e^{3} \right)$ in lowest terms, where $\ord_{2} \left( b_{1} \right)
\equiv 0 \bmod 2$ and $a=b_{1}b_{2}$ with $\gcd(M,N)=\gcd \left( e,N \right)=1$.

Suppose $b_{1}$ and $N$ are both even. Since $P \in E_{a}(\bQ)$, we can
write $N^{2} = b_{1}M^{4} + b_{2}e^{4}$. As $\gcd \left( e,N \right)=1$ and
since $\ord_{2} \left( b_{1} \right) \geq 2$, it follows that $p^{2}|b_{2}$.
Therefore $p^{4}$ divides $a$, which contradicts the assumption that $a$ is
fourth-power-free. Hence $b_{1}$ and $N$ cannot both be even.

We apply \eqref{eq:x2p-b2} to complete the proof for $\ord_{2}(x(P)) \equiv 0 \bmod 2$.

If $N$ is odd, then $2b_{1}M^{4}-N^{2}$ is odd and so
$\ord_{2} \left( B_{2} \right) \geq \ord_{2} \left( e^{2} \right)+2
= \ord_{2} \left( B_{1} \right)+2$.

If $N$ is even, then $b_{1}M^{4}$ is odd, since we saw that $b_{1}$ must be odd
and $\gcd(M,N)=1$. Hence, $2^{2}|| \left( 2b_{1}M^{4}-N^{2} \right)^{2}$, but
$2^{4}|4M^{2}N^{2}$. So, in this case too,
$\ord_{2} \left( B_{2} \right) \geq \ord_{2} \left( e^{2} \right)+2
= \ord_{2} \left( B_{1} \right)+2$.

Finally, if $\ord_{2}(x(P)) \not\equiv 0 \bmod 2$, then $\ord_{2}(x(P))
= \ord_{2} \left( b_{1} \right)$ is an odd positive integer. Because $a=b_{1}b_{2}
\equiv 4, 52 \bmod 64$, we must have $\ord_{2} \left( b_{1} \right)=1$,
as required.
\end{proof}

\section{Proof of Theorem~\ref{thm:lang2}}

We compute the canonical height by summing local heights.

Writing $x(2P)=\alpha^{2}/\delta^{2}$ as a fraction in lowest terms ($x(2P)$
is so expressible via Lemma~\ref{lem:2.4}(a)) and taking the sum of
\eqref{eq:padic-2p-hgt} and \eqref{eq:2adic-2p-hgt} over all primes gives
the exact formula
\begin{equation}
\label{eq:nonarch-sum}
\sum_{v \neq \infty} \widehat{\lambda}_{v}(2P)
= \log |\delta| + \frac{1}{12} \log \left| \Delta \left( E_{a} \right) \right|
- \left\{
	\begin{array}{ll}
		\frac{1}{2} \log(2) & \mbox{if $a \equiv 4, 52 \bmod 64$ and $\ord_{2}(x(2P))>0$} \\
		0 & \mbox{otherwise.}
	\end{array}
\right.
\end{equation}

\noindent
$\bullet$ $a<0$

Adding \eqref{eq:nonarch-sum} to the lower bound in Lemma~\ref{lem:arch-a-neg}(a)
for $\widehat{\lambda}_{\infty}(2P)$, we obtain
$$
\widehat{h}(2P)
> \frac{1}{4} \log \left| x^{2}(2P)-a \right| + \log |\delta|
- \left\{
	\begin{array}{ll}
		\frac{1}{2} \log(2) & \mbox{if $a \equiv 4, 52 \bmod 64$ and $\ord_{2}(x(2P))>0$} \\
		0 & \mbox{otherwise.}
	\end{array}
\right.
$$

Since $2P \in E^{0}(\bR)$, $x(2P) \geq \sqrt{|a|}$ and therefore
$x^{2}(2P)-a \geq |2a|$. Thus
$$
\widehat{h}(2P)
> \frac{1}{4} \log |2a| + \log |\delta|
- \left\{
	\begin{array}{ll}
		\frac{1}{2} \log(2) & \mbox{if $a \equiv 4, 52 \bmod 64$ and $\ord_{2}(x(2P))>0$} \\
		0 & \mbox{otherwise.}
	\end{array}
\right.
$$

Theorem~\ref{thm:lang2} follows in this case since $\widehat{h}(2P)=4\widehat{h}(P)$
and using Lemma~\ref{lem:2.4}(b) to provide a lower bound for $\log| \delta|$.

\noindent
$\bullet$ $a>0$

Here we proceed similarly. Adding \eqref{eq:nonarch-sum} to the lower bound
in Lemma~\ref{lem:arch-a-pos}(a), using Lemma~\ref{lem:2.4}(b) to provide a
lower bound for $\log| \delta|$, Theorem~\ref{thm:lang2} for $a>0$ follows
since $\widehat{h}(2P)=4\widehat{h}(P)$.

\section{Proof of Theorem~\ref{thm:hgt-diff}}

From Lemmas~\ref{lem:nonarch-v-odd}(c) and \ref{lem:nonarch-v-2}(c),
\begin{eqnarray}
\label{eq:non-arch}
0 & \leq &
\sum_{v \neq \infty} \left(
\frac{1}{2} \log \max \left\{ 1, |x(P)|_{v} \right\}
- \frac{1}{12} \log \left| \Delta \left( E_{a} \right) \right|_{v}
- \widehat{\lambda}_{v}(P) \right) \\
& \leq & \frac{1}{4}\log|a| + \frac{3}{8} \log(2), \nonumber
\end{eqnarray}
with the upper bound achieved when $a \equiv 4, 52 \bmod 64$ and $\ord_{2}(x(P))=1$,
by observing that $(7/8)\log(2) = (1/4)\log \left| 2^{\ord_{2}(a)} \right| + (3/8)\log(2)$
for such $a$.

If $a \leq -2$ and $P \in E_{a}^{0}(\bQ)$, then from Lemma~\ref{lem:arch-a-neg}(b)
and \eqref{eq:non-arch}
\begin{equation}
\label{eq:a-neg-p-pos}
-\frac{1}{3}\log(2)
< \frac{1}{2}h(P) - \widehat{h}(P)
< \frac{1}{4}\log|a| + \frac{3}{8} \log(2).
\end{equation}

Next suppose that $a \leq -2$ and $P \not\in E_{a}^{0}(\bQ)$, then from
Lemma~\ref{lem:arch-a-neg}(c) and \eqref{eq:non-arch}
\begin{equation}
\label{eq:a-neg-p-neg}
-\frac{1}{4} \log |a| - \frac{1}{2\sqrt{|a|}}
< \frac{1}{2}h(P) - \widehat{h}(P)
< \frac{1}{4}\log|a| + \frac{1}{8}\log(2).
\end{equation}

Now suppose that $a \geq 2$, then from Lemma~\ref{lem:arch-a-pos}(b) and
\eqref{eq:non-arch}
\begin{equation}
\label{eq:a-pos}
-\frac{1}{4}\log(a) - \frac{1}{2\sqrt{a}}
< \frac{1}{2}h(P) - \widehat{h}(P)
< \frac{1}{4}\log|a| + \frac{3}{8} \log(2).
\end{equation}

The upper bound in the first inequality in our Theorem is established.
The lower bound in the first inequality in our Theorem comes from combining the
lower bounds in \eqref{eq:a-neg-p-pos}, \eqref{eq:a-neg-p-neg} and
\eqref{eq:a-pos}, and noting that
$-(1/4)\log |a| -1/(2\sqrt{|a|} < -(1/3)\log(2)$ for all $a \leq -2$.

The lower bound in the second inequality in Theorem~\ref{thm:hgt-diff} comes
from using Lemma~\ref{lem:arch-a-neg}(d) and Lemma~\ref{lem:arch-a-pos}(c) above
rather than Lemma~\ref{lem:arch-a-neg}(c) and Lemma~\ref{lem:arch-a-pos}(b).
Once again, we note that
$-(1/4)\log |a| - 0.16 < -(1/3)\log(2)$ for all $a \leq -2$.

For $a=-1,1$ and $2$, $E_{a}(\bQ)$ consists only of torsion points,
which we consider next for all $a$.

From Proposition~X.6.1(a) of \cite{Silv2}, the torsion group of $E_{a}(\bQ)$ is
isomorphic to:\\
--$\bZ/4\bZ$, if $a=4$ (in this case the torsion points are $(0,0)$, $(2,\pm 4)$, $O$) \\
--$\bZ/2\bZ \times \bZ/2\bZ$, if $-a$ is a perfect square (in this case the torsion
points are $(\pm \sqrt{-a},0)$, $(0,0)$ and $O$) \\
--$\bZ/2\bZ$, otherwise (here the torsion points are $(0,0)$ and $O$).

In each of these cases, $(1/2)h(P)-\widehat{h}(P)=(1/4)\log|a|$ or $0$ and so
our Theorem holds for the torsion points too.

\section{Sharpness of Results}
\label{sect:sharp}

In the introduction, we stated that Theorems~\ref{thm:lang2} and \ref{thm:hgt-diff}
are best possible. We justify these statements here by constructing examples
demonstrating this.

\subsection{Lower Bounds}

Our non-archimedean results are exact, so any gap between the actual height
of points and our lower bounds in Theorem~\ref{thm:lang2} must arise from
the archimedean local height.

\noindent
$\bullet$ $a>0$

For $a>0$, if $x(2P)/\sqrt{a} \rightarrow 0$ as $\sqrt{a}$ grows, then the
left-hand side of \eqref{eq:archim-hgt-a-pos2} minus its right-hand side
approaches $0$. So we shall choose $x(2P)$ small and fixed and then determine
$a$ such that $2P \in E_{a}(\bQ)$.

We find $a$ in two steps. First, we determine for which values of $a$,
there exists a point $Q \in E_{a}(\bQ)$ with  $x(Q)=1/16$ for
$a \equiv 1,5,7,9,13,15 \bmod 16$, $x(Q)=1/4$ for $a \equiv 2,3,6,8,10,11,12,14 \bmod 16$
or $a \equiv 20,36 \bmod 64$, $x(Q)=4$ for $a \equiv 52 \bmod 64$ and
$x(Q)=16$ for $a \equiv 4 \bmod 64$. Such points $Q$ are suggested by
Lemma~\ref{lem:2.4}(b). For each congruence class, we obtain a pair of
quadratic polynomials such that $a$ must take the value of one of these
polynomials.

Second, we determine for which values of $a$ (restricted to values of the
aforementioned quadratic polynomials), such $Q$ are, in fact, $2P$.
In this way, we obtain the expressions in Table~\ref{table:lang-egs-a-pos}.
As $a_{1} \rightarrow +\infty$, the difference between the height of the
point $P$ on $E_{a}$ and the lower bound in Theorem~\ref{thm:lang2} approaches
$0$.

\noindent
$\bullet$ $a<0$

For $a<0$, the analysis is somewhat more complicated since $x(2P)$ must be close
to $\sqrt{|a|}$, which is not always rational. But we can proceed in a similar way.

For $a \equiv 2,3,6,8,10,11,12,14 \bmod 16$ or $a \equiv 20,36 \bmod 64$, let
$c$ be an odd positive integer, put $a=-\left( c^{4}-1 \right)/16$ and
$Q=\left( c^{2}/4, c/8 \right)$. Observe that $Q \in E_{a}(\bQ)$.

Now we want to find values of $c$ such that $Q=2P$ for some $P \in E_{a}(\bQ)$.
If $Q=2P$, then $x(2P)=\left( 16x(P)^{2}+c^{4}-1 \right)^{2}
/\left( 64x(P)\left( 16x(P)^{2}-c^{4}+1 \right) \right)=c^{2}/4$. This simplifies to
$$
\left( 16x(P)^{2}-8\left(c^{2}+1 \right)x(P) -c^{4}+1 \right)
\left( 16x(P)^{2}-8\left(c^{2}-1 \right)x(P) -c^{4}+1 \right)=0.
$$

One of these quadratic polynomials will have rational roots if and only if
$c^{2}-2z^{2}=\pm 1$ for some $z \in \bZ$. Such $c$ are members of the recurrence
sequence defined by $c_{n}=2c_{n-1}+c_{n-2}$ where $c_{0}=c_{1}=1$.

Lastly, we must ensure that the resulting values of $a$ are in the correct
congruence classes (as two examples, we require $c \equiv 15,49 \bmod 256$ for
$a \equiv 36,20 \bmod 64$, respectively). This imposes the restriction that the
index $n$ on the recurrence sequences must lie in certain congruence classes.

For $a \equiv 1,5,7,9,13,15 \bmod 16$, we proceed in the same way. Let
$a=-\left(c^{4}-1 \right)/256$, and $Q=\left( c^{2}/16, c/64 \right)$.
Here the polynomial used to ensure that $Q=2P$ is
$$
\left( 256x(P)^{2}-32\left(c^{2}+1 \right)x(P) -c^{4}+1 \right)
\left( 256x(P)^{2}-32\left(c^{2}-1 \right)x(P) -c^{4}+1 \right)
$$
and this implies that $c$ is a member of the same recurrence sequence as above.
Again, we must restrict the index $n$ on the recurrence sequences to certain
congruence classes. These congruences have larger moduli here since we also
require $c \equiv 1 \bmod 64$ so that our expression for $a$ is an integer.

For $a \equiv 4,52 \bmod 64$, we proceed similarly. Let $a=-\left(d^{4}-4 \right)$,
and $Q=\left( d^{2}, 2d \right)$. If $d \equiv 0 \bmod 4$, then $a \equiv 4 \bmod 64$;
while if $d \equiv 2 \bmod 4$, then $a \equiv 52 \bmod 64$. To ensure that $Q=2P$,
we use the polynomial
$$
\left( x(P)^{2}+\left(-2d^{2}-4 \right)x(P) -d^{4}+4 \right)
\left( x(P)^{2}+\left(-2d^{2}+4 \right)x(P) -d^{4}+4 \right).
$$

Here we require $2d^{2}-z^{2}=\pm 4$, for some $z \in \bZ$. Hence $d$ are
members of the recurrence sequence defined by $d_{n}=2d_{n-1}+d_{n-2}$ where
$d_{0}=0$ and $d_{1}=1$.

We find that $d_{4n} \equiv 0 \bmod 4$ and $d_{4n+2} \equiv 2 \bmod 4$.

\begin{table}[ht]
\centering
\begin{tabular}{|r|r|r|} \hline
                       &   $a$                                                                                     & $x(P)$  \\ \hline
$a \equiv  1 \bmod 16$ & $ \left(  16a_{1}+ 1 \right) \left( 256a_{1}+ 17 \right) \left( 512a_{1}+ 33 \right)^{2}$ & $ \left( 16a_{1}+ 1 \right) \left( 512a_{1}+ 33 \right)$ \\
$a \equiv  2 \bmod 16$ & $2\left(  16a_{1}+ 7 \right) \left(  32a_{1}+ 15 \right) \left(  64a_{1}+ 29 \right)^{2}$ & $ \left( 32a_{1}+15 \right) \left(  64a_{1}+ 29 \right)$ \\
$a \equiv  3 \bmod 16$ & $ \left(  32a_{1}+13 \right) \left(  32a_{1}+ 15 \right) \left(  16a_{1}+  7 \right)^{2}$ & $ \left( 16a_{1}+ 7 \right) \left(  32a_{1}+ 15 \right)$ \\
$a \equiv  5 \bmod 16$ & $ \left(  16a_{1}+ 5 \right) \left( 256a_{1}+ 81 \right) \left( 512a_{1}+161 \right)^{2}$ & $4\left( 16a_{1}+ 5 \right) \left( 512a_{1}+161 \right)$ \\
$a \equiv  6 \bmod 16$ & $2\left(  16a_{1}+ 7 \right) \left(  32a_{1}+ 13 \right) \left(  64a_{1}+ 27 \right)^{2}$ & $ \left( 32a_{1}+13 \right) \left(  64a_{1}+ 27 \right)$ \\
$a \equiv  7 \bmod 16$ & $ \left(  64a_{1}+23 \right) \left(  64a_{1}+ 25 \right) \left(   8a_{1}+  3 \right)^{2}$ & $ \left( 64a_{1}+25 \right) \left(   8a_{1}+  3 \right)$ \\
$a \equiv  8 \bmod 16$ & $8\left(   2a_{1}+ 1 \right) \left(  16a_{1}+  7 \right) \left(  32a_{1}+ 15 \right)^{2}$ & $ \left( 16a_{1}+ 7 \right) \left(  32a_{1}+ 15 \right)$ \\
$a \equiv  9 \bmod 16$ & $ \left(  16a_{1}+ 7 \right) \left( 256a_{1}+111 \right) \left( 512a_{1}+223 \right)^{2}$ & $4\left( 16a_{1}+ 7 \right) \left( 512a_{1}+223 \right)$ \\
$a \equiv 10 \bmod 16$ & $2\left(   8a_{1}+ 3 \right) \left(  16a_{1}+  7 \right) \left(  32a_{1}+ 13 \right)^{2}$ & $ \left( 16a_{1}+ 7 \right) \left(  32a_{1}+ 13 \right)$ \\
$a \equiv 11 \bmod 16$ & $ \left(  16a_{1}+ 5 \right) \left(  16a_{1}+  7 \right) \left(   8a_{1}+  3 \right)^{2}$ & $4\left( 16a_{1}+ 5 \right) \left(   8a_{1}+  3 \right)$ \\
$a \equiv 12 \bmod 16$ & $4\left(   4a_{1}+ 1 \right) \left(  16a_{1}+  3 \right) \left(  32a_{1}+  7 \right)^{2}$ & $ \left( 16a_{1}+ 3 \right) \left(  32a_{1}+  7 \right)$ \\
$a \equiv 13 \bmod 16$ & $ \left(  16a_{1}+ 3 \right) \left( 256a_{1}+ 47 \right) \left( 512a_{1}+ 95 \right)^{2}$ & $4\left( 16a_{1}+ 3 \right) \left( 512a_{1}+ 95 \right)$ \\
$a \equiv 14 \bmod 16$ & $2\left(   8a_{1}+ 3 \right) \left(  16a_{1}+  5 \right) \left(  32a_{1}+ 11 \right)^{2}$ & $ \left( 16a_{1}+ 5 \right) \left(  32a_{1}+ 11 \right)$ \\
$a \equiv 15 \bmod 16$ & $ \left( 128a_{1}+55 \right) \left( 128a_{1}+ 57 \right) \left(  16a_{1}+  7 \right)^{2}$ & $ \left( 16a_{1}+ 7 \right) \left( 128a_{1}+ 55 \right)$ \\
$a \equiv  4 \bmod 64$ & $4\left(   2a_{1}- 1 \right) \left(   2a_{1}+  7 \right) \left(   2a_{1}+  3 \right)^{2}$ & $2\left(  2a_{1}+ 3 \right) \left(   2a_{1}+  7 \right)$ \\
$a \equiv 20 \bmod 64$ & $4\left(  16a_{1}- 1 \right) \left(  64a_{1}-  5 \right) \left( 128a_{1}-  9 \right)^{2}$ & $4\left( 16a_{1}- 1 \right) \left( 128a_{1}-  9 \right)$ \\
$a \equiv 36 \bmod 64$ & $4\left(  16a_{1}-11 \right) \left(  64a_{1}- 43 \right) \left( 128a_{1}- 87 \right)^{2}$ & $ \left( 64a_{1}-43 \right) \left( 128a_{1}- 87 \right)$ \\
$a \equiv 52 \bmod 64$ & $4\left(   8a_{1}+ 1 \right) \left(   8a_{1}+  5 \right) \left(   8a_{1}+  3 \right)^{2}$ & $2\left(  8a_{1}+ 1 \right) \left(   8a_{1}+  3 \right)$ \\ \hline
\end{tabular}
\caption{$\widehat{h}(P)$ near bounds for $a>0$}
\label{table:lang-egs-a-pos}
\end{table}

\begin{table}[ht]
\centering
\begin{tabular}{|r|r|r|} \hline
                       & $a$                                      & $x(2P)$               \\ \hline
$a \equiv  1 \bmod 16$ & $-\left( c_{512n+161}^{4}-1 \right)/256$ & $c_{512n+161}^{2}/16$ \\
$a \equiv  2 \bmod 16$ & $-\left( c_{32n+13}^{4}-1 \right)/16$    & $c_{32n+13}^{2}/4$    \\
$a \equiv  3 \bmod 16$ & $-\left( c_{16n+6}^{4}-1 \right)/16$     & $c_{16n+6}^{2}/4$     \\
$a \equiv  5 \bmod 16$ & $-\left( c_{512n+289}^{4}-1 \right)/256$ & $c_{512n+289}^{2}/16$ \\
$a \equiv  6 \bmod 16$ & $-\left( c_{32n+11}^{4}-1 \right)/16$    & $c_{32n+11}^{2}/4$    \\
$a \equiv  7 \bmod 16$ & $-\left( c_{64n+8}^{4}-1 \right)/256$    & $c_{64n+8}^{2}/16$    \\
$a \equiv  8 \bmod 16$ & $-\left( c_{32n+15}^{4}-1 \right)/16$    & $c_{32n+15}^{2}/4$    \\
$a \equiv  9 \bmod 16$ & $-\left( c_{512n+417}^{4}-1 \right)/256$ & $c_{512n+417}^{2}/16$ \\
$a \equiv 10 \bmod 16$ & $-\left( c_{32n+3}^{4}-1 \right)/16$     & $c_{32n+3}^{2}/4$     \\
$a \equiv 11 \bmod 16$ & $-\left( c_{16n+2}^{4}-1 \right)/16$     & $c_{16n+2}^{2}/4$     \\
$a \equiv 12 \bmod 16$ & $-\left( c_{16n+4}^{4}-1 \right)/16$     & $c_{16n+4}^{2}/4$     \\
$a \equiv 13 \bmod 16$ & $-\left( c_{512n+33}^{4}-1 \right)/256$  & $c_{512n+33}^{2}/16$  \\
$a \equiv 14 \bmod 16$ & $-\left( c_{32n+5}^{4}-1 \right)/16$     & $c_{32n+5}^{2}/4$     \\
$a \equiv 15 \bmod 16$ & $-\left( c_{64n+24}^{4}-1 \right)/256$   & $c_{64n+24}^{2}/16$   \\
$a \equiv  4 \bmod 64$ & $-\left( d_{4n}^{4}-4 \right)$           & $d_{4n}^{2}$          \\
$a \equiv 20 \bmod 64$ & $-\left( c_{128n+41}^{4}-1 \right)/16$   & $c_{128n+41}^{2}/4$   \\
$a \equiv 36 \bmod 64$ & $-\left( c_{128n+55}^{4}-1 \right)/16$   & $c_{128n+55}^{2}/4$   \\
$a \equiv 52 \bmod 64$ & $-\left( d_{4n+2}^{4}-4 \right)$         & $d_{4n+2}^{2}$        \\ \hline
\end{tabular}
\caption{$\widehat{h}(P)$ near bounds for $a<0$}
\end{table}

\subsection{Difference of Heights}

Silverman (see Example~2.2 of \cite{Silv5}) shows that the coefficients of the
$\log |a|$ terms are best possible.

As for the constants, $-0.16$ cannot be replaced by anything greater than
$-0.14913\ldots$ by considering the point $1158(1,2)$ on $y^{2}=x^{3}+3x$
(note that $x(1158(1,2))=0.999402\ldots$). Taking the archimedean height function
evaluated at $x=1$ for $a=3$, we see that $-0.14922$ is the smallest possible
constant ($-0.1310\ldots$ appears to be the worst value for $a<0$, by
considering the point $(-1,1)$ on $y^{2}=x^{3}-2x$).

However, the correct constant there is $0$, so the constant in our lower bound
is best possible. Taking $a=-4a_{1}^{2}-4a_{1}-2$ and $P= \left(-1,2a_{1}+1 \right)$
provides examples for $a<0$, while $a=4a_{1}^{2}+4a_{1}$ and $P=\left(1,2a_{1}+1 \right)$
provides examples for $a>0$.

The constant $(3/8)\log(2)=0.25993019\ldots$ is best possible as is demonstrated
by considering $a=32a_{1}^{2}+32a_{1}+4$, $P=(a/2,4(2a_{1}+1)(8a_{1}^{2}+8a_{1}+1))$
and letting $a_{1}$ be an arbitrarily large integer. For example, with $a_{1}=10^{6}$,
the constant is $0.25993011\ldots$.

\bibliographystyle{amsplain}

\begin{thebibliography}{20}
\bibitem{BG}
E. Bombieri, W. Gubler,
\emph{Heights in Diophantine Geometry},
Cambridge University Press, 2006.

\bibitem{Bremner}
A. Bremner, J.H. Silverman, N. Tzanakis,
\emph{Integral points in arithmetic progression on $y^{2}=x\left( x^{2}-n^{2} \right)$},
J. Number Theory 80(2) (2000), 187--208. 

\bibitem{Cremona}
J. E. Cremona, M. Prickett, S. Siksek,
\emph{Height difference bounds for elliptic curves over number fields},
J. Number Theory 116 (2006), 42--68.

\bibitem{David}
S. David,
\emph{Points de petite hauteur sur les courbes elliptiques},
J. Number Theory 64 (1997) 104--129.

\bibitem{Everest1}
G. Everest, G. Mclaren, T. Ward,
\emph{Primitive divisors of elliptic divisibility sequenses},
J. Number Theory 118 (2006), 71--89.

\bibitem{EIS}
G. Everest, P. Ingram, S. Stevens,
\emph{Primitive divisors on twists of the Fermat cubic},
LMS Journal of Computation and Mathematics 12 (2009), 54--81.

\bibitem{HS}
M. Hindry, J. H. Silverman,
\emph{The canonical height and integral points on elliptic curves},
Invent. Math. 93 (1988), 419--450.

\bibitem{Lang}
S. Lang,
\emph{Elliptic Curves: Diophantine analysis},
Grundlehren der Mathematischen Wissenschaften 231, Springer-Verlag, Berlin, 1978.

\bibitem{Masser}
D. Masser,
\emph{Counting points of small height on elliptic curves},
Bull. Soc. Math. France 117 (1989), 247--265.

\bibitem{Petsche}
C. Petsche,
\emph{Small rational points on elliptic curves over number fields},
New York Journal of Math. 12 (2006), 257--268.

\bibitem{S-Z}
S. Schmitt, H. G. Zimmer,
\emph{Elliptic Curves: A Computational Approach},
de Gruyter Studies in Mathematics 31, Walter de Gruyter Inc, 2004.

\bibitem{Silv1}
J.~H. Silverman,
\emph{Lower bound for the canonical height on elliptic curves},
Duke Math. J. 48 (1981), 633--648.

\bibitem{Silv2}
J.~H. Silverman,
\emph{The Arithmetic of Elliptic Curves},
Graduate Texts in Math. 106, Springer-Verlag, New York, 1986.

\bibitem{Silv3}
J.~H. Silverman.
\emph{Computing heights on elliptic curves},
Math. Comp. 51 (1988), 339--358.

\bibitem{Silv5}
J.~H. Silverman,
\emph{The difference between the Weil height and the canonical height on elliptic curves},
Math. Comp. 55 (1990), 723--743.

\bibitem{Silv6}
J.~H. Silverman,
\emph{Advanced Topics in the Arithmetic of Elliptic Curves},
Graduate Texts in Math. 151, Springer-Verlag, New York, 1994.

\bibitem{Silv-Tate}
J.~H. Silverman, J. Tate,
\emph{Rational Points on Elliptic Curves},
Undergraduate Texts in Math., Springer-Verlag, New York, 1992.

\bibitem{Tate}
J. Tate,
Letter to J.-P. Serre, 1 Oct 1979 (see \verb+http://arxiv.org/abs/1207.5765+).

\bibitem{VY1}
P. Voutier and M. Yabuta,
\emph{Primitive divisors of certain elliptic divisibility sequences},
Acta Arith. 151 (2012), 165--190.
\end{thebibliography}

\end{document}